\documentclass[12pt]{amsart}

\pdfoutput=1

\usepackage{amsmath, amsthm, amssymb}
\usepackage{url} 
\usepackage{graphicx} 

\theoremstyle{plain}
\newtheorem{Theorem}{Theorem}[section]
\newtheorem{Lemma}[Theorem]{Lemma}

\newtheorem{Proposition}[Theorem]{Proposition}
\newtheorem{Conjecture}[Theorem]{Conjecture}
\theoremstyle{remark}
\newtheorem{Remark}[Theorem]{Remark}

\theoremstyle{definition}
\newtheorem{Example}[Theorem]{Example}

\newtheorem{Question}[Theorem]{Question}

\def\Hc{\mathcal{H}}
\def\Pc{\mathcal{P}}
\def\ZZ{\mathbb{Z}}
\def\RR{\mathbb{R}}
\def\CC{\mathbb{C}}
\def\eb{\mathbf{e}}
\def\Vb{\mathbf{V}}
\def\a{\alpha}
\def\d{\delta}

\let\phi=\varphi

\let\union=\cup
\let\Union=\bigcup

\def\th{-th }

\def\cocoa{{\hbox{\rm C\kern-.13em o\kern-.07em C\kern-.13em o\kern-.15em A}}}

\DeclareMathOperator{\con}{conv}
\DeclareMathOperator{\vol}{vol}

\renewcommand{\Re}{\operatorname{Re}}

\newcommand{\numberof}[1]{\left|#1\right|}
\newcommand{\floorof}[1]{\left\lfloor#1\right\rfloor}
\newcommand{\realpartof}[1]{\Re(#1)}

\begin{document}

\title{Roots of Ehrhart polynomials arising from graphs}
\author{Tetsushi Matsui, Akihiro Higashitani, Yuuki Nagazawa, Hidefumi Ohsugi and Takayuki Hibi}
\date{}

\thanks{
{\bf 2010 Mathematics Subject Classification:}
Primary 52C07; Secondary 52B20, 12D10. \\
\, \, \, {\bf Keywords:}
Ehrhart polynomial, edge polytope, Fano polytope, smooth polytope.
}
\address{Tetsushi Matsui,
Department of Pure and Applied Mathematics,
Graduate School of Information Science and Technology,
Osaka University,
Toyonaka, Osaka 560-0043, Japan.
/ JST CREST.\\
(Currently the author is at National Institute for Informatics,
Chiyoda-ku, Tokyo 101-8430, Japan.)
}
\email{t-matsui@cr.math.sci.osaka-u.ac.jp}
\address{Akihiro Higashitani,
Department of Pure and Applied Mathematics,
Graduate School of Information Science and Technology,
Osaka University,
Toyonaka, Osaka 560-0043, Japan}
\email{sm5037ha@ecs.cmc.osaka-u.ac.jp}
\address{Yuuki Nagazawa,
Department of Pure and Applied Mathematics,
Graduate School of Information Science and Technology,
Osaka University,
Toyonaka, Osaka 560-0043, Japan}
\email{sm5032ny@ecs.cmc.osaka-u.ac.jp}
\address{Hidefumi Ohsugi,
Department of Mathematics,
College of Science,
Rikkyo University,
Toshima-ku, Tokyo 171-8501, Japan.
/ JST CREST.}
\email{ohsugi@rikkyo.ac.jp}
\address{Takayuki Hibi,
Department of Pure and Applied Mathematics,
Graduate School of Information Science and Technology,
Osaka University,
Toyonaka, Osaka 560-0043, Japan.
/ JST CREST.}
\email{hibi@math.sci.osaka-u.ac.jp}

\begin{abstract}
Several polytopes arise from finite graphs. For edge and symmetric edge
polytopes, in particular, exhaustive computation of the Ehrhart polynomials
not merely supports the conjecture of Beck {\it et al.}\ that
all roots $\a$ of Ehrhart polynomials of polytopes of dimension $D$ satisfy
$-D \le \realpartof{\a} \le D - 1$,
but also reveals some interesting phenomena for each type of polytope.
Here we present two new conjectures:
(1) the roots of the Ehrhart polynomial of an edge polytope for a complete
multipartite graph of order $d$ lie in the circle 
$|z+\tfrac{d}{4}| \le \tfrac{d}{4}$ or are negative integers, and
(2) a Gorenstein Fano polytope of dimension $D$
has the roots of its Ehrhart polynomial in the narrower strip
$-\tfrac{D}{2} \leq \realpartof{\a} \leq \tfrac{D}{2}-1$.
Some rigorous results to support them are obtained as well as
for the original conjecture.
The root distribution of Ehrhart polynomials of each type of
polytope is plotted in figures.
\end{abstract}

\maketitle                        

\section*{Introduction}

The root distribution of Ehrhart polynomials is 
one of the current topics
on computational commutative algebra.
It is well-known that the coefficients of an Ehrhart polynomial
reflect combinatorial and geometric properties such as
the volume of the polytope in the leading coefficient,
gathered information about its faces in the second coefficient, etc.
The roots of an Ehrhart polynomial should also reflect
properties of a polytope that are hard to elicit just from the coefficients.
Among the many papers on the topic,
including~\cite{BHW2007},~\cite{Bra2008},~\cite{BD2006},~\cite{HSW2005}
and~\cite{Pfe2010},
Beck {\it et al.\/}~\cite{BDDPS2005} conjecture that:
\begin{Conjecture}\label{conj:dstrip}
  All roots $\a$ of Ehrhart polynomials of lattice $D$-polytopes
  satisfy $-D \le \realpartof{\a} \le D - 1$.
\end{Conjecture}
Compared with the norm bound, which is \(O(D^2)\) in general~\cite{Bra2008},
the strip in the conjecture puts a tight restriction on the
distribution of roots for any Ehrhart polynomial.

This paper investigates the roots of Ehrhart polynomials of polytopes
arising from graphs,
namely, edge polytopes and symmetric edge polytopes.
The results obtained not merely support Conjecture~\ref{conj:dstrip},
but also reveal some interesting phenomena.
Regarding the scope of the paper, note
that both kinds of polytopes are ``small'' in a sense: That is,
each edge polytope from a graph without loops is contained in a unit hypercube,
and one from a graph with loops, in twice a unit hypercube;
whereas each symmetric edge polytope is contained in twice a unit hypercube.

In Section~\ref{sec:simple},
the distribution of roots of Ehrhart polynomials of edge polytopes 
is computed, and as a special case, that of complete multipartite graphs
is studied.
We observed from exhaustive computation that all roots have a negative
real part and they are in the range of Conjecture~\ref{conj:dstrip}.
Moreover, for complete multipartite graphs of order $d$,
the roots lie in the circle
$|z+\tfrac{d}{4}|\le\tfrac{d}{4}$ or are negative integers greater than $-(d-1)$.
And we conjecture its validity beyond the computed range of $d$
(Conjecture~\ref{conj:circle}).

Simple edge polytopes constructed from graphs with possible loops are
studied in Section~\ref{sec:loops}.
Roots of the Ehrhart polynomials are determined in some cases.
Let $G$ be a graph of order $d$ with loops and
$G'$ its subgraph of order $p$ induced by
vertices without a loop attached.
Then, Theorem~\ref{root1} proves that the real roots are
in the interval $[-(d-2),0)$, especially all integers in
$\{-(d-p),\ldots,-1\}$ are roots of the polynomial;
Theorem~\ref{root2} determines that if $d -2p +2 \geq 0$,
there are $p-1$ real non-integer roots each of which is unique
in one of ranges $(-k, -k+1)$ for $k=1,\ldots,p-1$;
and Theorem~\ref{root3} proves that if $d > p \geq 2$,
all the integers $-\floorof{\frac{d-1}{2}},\ldots,-1$ are
roots of the polynomial.
We observed that all roots have a negative real part and
are in the range of Conjecture~\ref{conj:dstrip}.

The symmetric edge polytopes in Section~\ref{sec:symmetric} are
Gorenstein Fano polytopes.
A unimodular equivalence condition for two symmetric edge polytopes is
also described in the language of graphs (Theorem~\ref{nonequivalent}).
The polytopes have Ehrhart polynomials with an interesting root distribution:
the roots are distributed symmetrically with respect to the vertical line
$\realpartof{z}=-\tfrac{1}{2}$.
We not only observe that all roots are in the range of
Conjecture~\ref{conj:dstrip}, but also conjecture that all roots
in $-\tfrac{D}{2} \leq \realpartof{\a} \leq \tfrac{D}{2}-1$
for Gorenstein Fano polytopes of dimension $D$
 (Conjecture~\ref{conj:narrowstrip}). 

\smallskip

Before starting the discussion, let us summarize the definitions of
edge polytopes, symmetric edge polytopes, etc.

Throughout this paper, graphs are always finite,
and so we usually omit the adjective ``finite.''
Let $G$ be a graph having no multiple edges on the vertex set
$V(G) = \{ 1, \ldots, d \}$
and the edge set
$E(G) = \{ e_1, \ldots, e_n \} \subset {V(G)}^2$.
Graphs may have loops in their edge sets unless explicitly excluded;
in which case the graphs are called {\em simple\/} graphs.
A {\em walk\/} of $G$ of length $q$ is a sequence 
$(e_{i_1}, e_{i_2}, \ldots, e_{i_q})$
of the edges of $G$, 
where 
$e_{i_k} = \{ u_k, u_{k+1} \}$
for $k = 1, \ldots, q$.
If, moreover, $u_{q+1} = u_{1}$ holds, then the walk is a {\em closed\/} walk.
Such a closed walk is called a {\em cycle\/} of length $q$ 
if $u_k \neq u_{k'}$ for all $1 \leq k < k' \leq q$.
In particular, a loop is a cycle of length 1.
Another notation, $(u_1,u_2,\ldots,u_q)$, will be also used for
the same cycle with $(\{u_1,u_2\},\{u_2,u_3\},\ldots,\{u_q,u_1\})$.
Two vertices $u$ and $v$ of $G$ are {\em connected\/}
if $u=v$ or there exists a walk $(e_{i_1}, e_{i_2}, \ldots, e_{i_q})$ of $G$
such that $e_{i_1} = \{u, v_1\}$ and $e_{i_q} = \{u_q, v\}$.
The connectedness is an equivalence relation and
the equivalence classes are called the {\em components\/} of $G$.
If $G$ itself is the only component, then $G$ is a {\em connected graph}.
For further information on graph theory,
we refer the reader to e.g.~\cite{Harary},~\cite{Wilson}

If $e = \{ i, j \}$ is an edge
of $G$ between $i \in V(G)$ and $j \in V(G)$, 
then we define $\rho(e) = \eb_i + \eb_j$.
Here, $\eb_i$ is the $i$\th unit coordinate 
vector of $\RR^d$.
In particular, for a loop 
$e = \{ i, i \}$ at $i \in V(G)$,
one has $\rho(e) = 2 \eb_i$.
The {\em edge polytope\/} of $G$  
is the convex polytope $\Pc_G$
$(\subset \RR^d)$, which is the convex hull of
the finite set 
$\{ \rho(e_1), \ldots, \rho(e_n) \}$.
The dimension of $\Pc_G$ equals to $d-2$ if the graph $G$ is a connected
bipartite graph,
or $d-1$, other connected graphs~\cite{OH1998}. 
The edge polytopes of complete multipartite graphs are studied in~\cite{OH2000}.
Note that if the graph $G$ is a complete graph,
the edge polytope $\Pc_G$ is also called the second hypersimplex
in~\cite[Section 9]{Sturmfels1995}.

Similarly, we define $\sigma(e) = \eb_i - \eb_j$
for an edge $e=\{ i, j \}$ of a simple graph $G$.
Then, the {\em symmetric edge polytope\/} of $G$
is the convex polytope $\Pc^{\pm}_G$
$(\subset \RR^d)$, which is the convex hull of
the finite set
$\{ \pm\sigma(e_1), \ldots, \pm\sigma(e_n) \}$.
Note that if $G$ is the complete graph $K_d$,
the symmetric edge polytope $\Pc^{\pm}_{K_d}$ coincides with
the root polytope of the lattice $A_d$ defined in~\cite{ABHPS2008}.

If $\Pc \subset \RR^N$ is an integral convex polytope,
then we define $i(\Pc, m)$ by
\[
i(\Pc,m) = \numberof{m \Pc \cap \ZZ^N}
.\]
We call $i(\Pc,m)$ the {\em Ehrhart polynomial\/} of $\Pc$
after Ehrhart, who succeeded in proving that $i(\Pc, m)$ is
a polynomial in $m$ of degree $\dim \Pc$ with $i(\Pc, 0) = 1$.
If $\vol(\Pc)$ is the normalized volume of $\Pc$, then
the leading coefficient of $i(\Pc,m)$ is $\tfrac{\vol(\Pc)}{(\dim \Pc)!}$.

An Ehrhart polynomial $i(\Pc, m)$ of $\Pc$ is related to a sequence of
integers called the {\em $\d$-vector\/},
$\d(\Pc) = (\d_0, \d_1, \ldots, \d_D)$, of $\Pc$ by
\[\sum_{m=0}^{\infty}i(\Pc, m)t^m = \frac{\sum_{j=0}^D\d_j t^j}{{(1-t)}^{D+1}}\]
where \(D\) is the degree of $i(\Pc, m)$.
We call the polynomial in the numerator on the right-hand side of the
equation above $\d_{\Pc}(t)$, the {\em $\d$-polynomial\/} of $\Pc$.
Note that the $\d$-vectors and $\d$-polynomials are referred to by other names
in the literature: e.g., in~\cite{Stanley1986},~\cite{Sta1993},
$h^*$-vector or $i$-Eulerian numbers are synonyms of $\d$-vector,
and $h^*$-polynomial or $i$-Eulerian polynomial, of $\d$-polynomial.
It follows from the definition that $\d_0 = 1$,
$\d_1 = \numberof{\Pc \cap \ZZ^N} - (D + 1)$, etc.
It is known that each $\d_i$ is nonnegative~\cite{Sta1980}.
If $\d_D \neq 0$, then
$\d_1 \leq \d_i$ for every $1 \leq i < D$~\cite{Hib1994}.
Though the roots of the polynomial are the focus of this paper,
the $\d$-vector is also a very important research subject.
For the detailed discussion on Ehrhart polynomials of convex polytopes,
we refer the reader to~\cite{Hibi1992}.

\section{Edge polytopes of simple graphs}
\label{sec:simple}

The aim in this section is to confirm Conjecture~\ref{conj:dstrip}
for the Ehrhart polynomials of edge polytopes constructed from
connected simple graphs, mainly by computational means.

\subsection{Exhaustive Computation for Small Graphs}
\label{sec:resul}

Let $\CC[X]$ denote the polynomial ring in one variable over
the field of complex numbers.  Given a polynomial 
$f = f(X) \in \CC[X]$, we write
$\Vb(f)$ for the set of roots of $f$, i.e.,
\[
\Vb(f) = \{ a \in \CC \ | \  f(a) = 0 \}.
\]

We computed the Ehrhart polynomial $i(\Pc_G, m)$ of each edge polytope
$\Pc_G$ for connected simple graphs $G$ of orders up to nine;
there are \(1, 2, \ldots, 261080\) connected simple graphs of orders
\(2, 3, \ldots, 9\)\footnote{
These numbers of such graphs are known;
see, e.g.,~\cite[Chapter 4]{HarPal} or {\tt A001349} of
the On-Line Encyclopedia of Integer Sequences.}.
Then, we solved each equation $i(\Pc_G, X)=0$ in the field of complex numbers.
For the readers interested in our method of computation,
see the small note in Appendix~\ref{sec:metho}.

\begin{figure}[htb!]
\centering%
\includegraphics[scale=0.4]{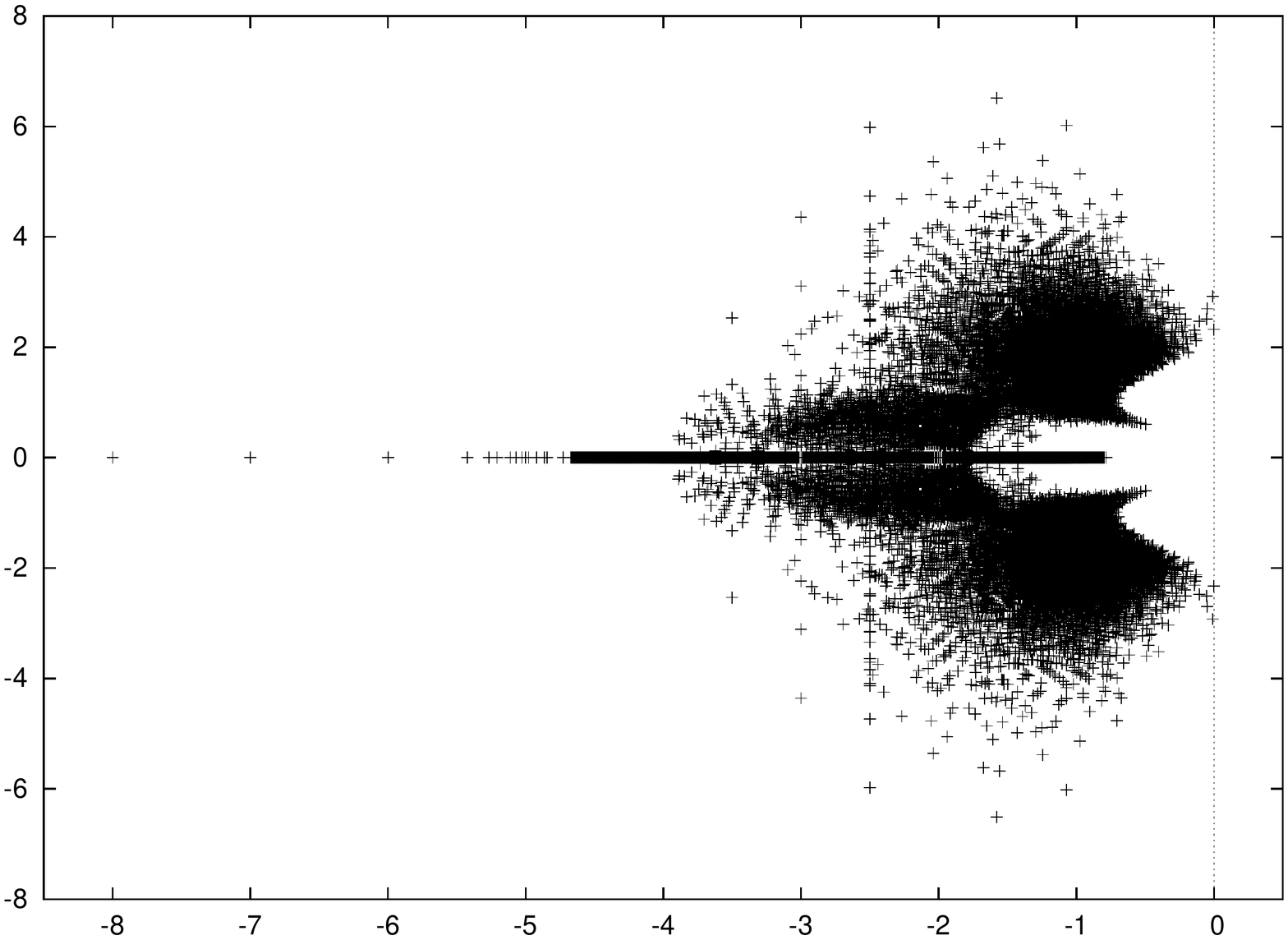}
\caption{$\Vb^{\text cs}_9$}
\label{fig:root9}
\end{figure}

Let $\Vb^{\text cs}_d$ denote $\Union \Vb(i(\Pc_G, m))$, 
where the union runs over all connected simple graphs $G$ of order $d$.
Figure~\ref{fig:root9} plots points of $\Vb^{\text cs}_9$,
as a representative of all results.
For all connected simple graphs of order $2$--$9$,
Conjecture~\ref{conj:dstrip} holds.

Since an edge polytope is a kind of $0/1$-polytope, 
the points in Figure~\ref{fig:root9} for $\Vb^{\text cs}_9$ 
are similar to those in Figure 6 of~\cite{BDDPS2005}.
However, the former has many more points, which form 
three clusters: one on the real axis, and other two being 
complex conjugates of each 
other and located nearer to the imaginary axis than the first cluster.
The interesting thing is that no roots appear in the right half plane
of the figure.
The closest points to the imaginary axis are approximately
$-0.583002 \pm 0.645775i \in \Vb^{\text cs}_7$,
$-0.213574 \pm 2.469065i \in \Vb^{\text cs}_8$, and
$-0.001610 \pm 2.324505i \in \Vb^{\text cs}_9$.
A polynomial with roots only in the left half plane
is called a {\em stable\/} polynomial.
This observation raises an open question:

\begin{Question}\label{q:leftplane}
  For any $d$ and any connected simple graph $G$ of order $d$,
  is $i(\Pc_G, m)$ always a stable polynomial?
\end{Question}

For a few infinite families of graphs, rigorous proofs are known:
Proposition~\ref{eg:complete} just below and Examples in the next subsection.

\begin{Proposition}
\label{eg:complete}
A root $\a$ of the Ehrhart polynomial $i(\Pc_{K_d}, m)$
of the complete graph $K_d$ satisfies
\begin{enumerate}
\item $\a \in \{-1, -2\}$ if $d=3$ or
\item $-\tfrac{d}{2} < \realpartof{\a} < 0$ if $d \ge 4$.
\end{enumerate}
\end{Proposition}
\begin{proof}
The Ehrhart polynomial $i(\Pc_{K_d}, m)$ of the complete graph $K_d$
is given in~\cite[Corollary~9.6]{Sturmfels1995}:
\[
i(\Pc_{K_d}, m)=\binom{d+2m-1}{d-1}-d\binom{m+d-2}{d-1} .
\]
In cases where $d=2$ or $3$, the Ehrhart polynomials are binomial coefficients,
since the edge polytopes are simplices.
Actually, they are:
\[
i(\Pc_{K_2}, m)=1\ \ \mbox{ and }\ \ 
i(\Pc_{K_3}, m)=\binom{m+2}{2}.
\]
Thus, there are no roots for $d=2$, whereas $\{-1, -2\}$ are the roots for $d=3$.

Hereafter, we assume $d \ge 4$.
It is easy to see that $\{-1,-2,\ldots,-\floorof{\tfrac{d-1}{2}}\}$
are included in $\Vb(i(\Pc_{K_d}, m))$.

We shall first prove that $\realpartof{\a} < 0$.
Let $q_d^{(1)}(m) = (2m+d-1)\cdots(2m+1)$ and $q_d^{(2)}(m)=d(m + d-2)\cdots m$.
Then for a complex number $z$, $i(\Pc_{K_d}, z)=0$ if and only if $q_d^{(1)}(z) = q_d^{(2)}(z)$, since $q_d^{(1)}(z) - q_d^{(2)}(z)$ is $(d-1)!\, i(\Pc_{K_d}, z)$.
Let us prove $|q_d^{(1)}(z)| > |q_d^{(2)}(z)|$ for any
complex number $z$ with a nonnegative real part
by mathematical induction on $d \geq 4$.

If $d=4$, 
\begin{eqnarray*}
  |q_4^{(1)}(z)| = |(2z + 3)(2z + 2)(2z + 1)|
  &=& |2z + 3||z + 1||4z + 2|\\
  &>& |z+2||z+1||4z|\, =\, |q_4^{(2)}(z)|
\end{eqnarray*}
holds for any complex number $z$ with $\realpartof{z} \ge 0$.

Assume for $d$ that $|q_d^{(1)}(z)| > |q_d^{(2)}(z)|$ is true for any
complex number $z$ with $\realpartof{z} \ge 0$.

Then, by
\begin{eqnarray*}
  |q_{d+1}^{(1)}(z)| &=& |2z + d||q_d^{(1)}(z)|\\
  |q_{d+1}^{(2)}(z)| &=& \dfrac{d+1}{d}|z + d - 1||q_d^{(2)}(z)|
\end{eqnarray*}
and
\begin{eqnarray*}
  |2dz + d^2| &>& |(d+1)z + d^2 - 1|
\end{eqnarray*}
from $2d > d+1$ and $d^2 > d^2 - 1$,
one can deduce
\begin{eqnarray*}
  d|q_{d+1}^{(1)}(z)| &=& |2dz + d^2||q_d^{(1)}(z)|\\
  &>& |(d+1)z + d^2 - 1||q_d^{(2)}(z)|\\
  &=& (d+1)|z + d - 1||q_d^{(2)}(z)|\\
  &=& d \dfrac{d+1}{d}|z + d - 1||q_d^{(2)}(z)|
  \, =\, d |q_{d+1}^{(2)}(z)|.
\end{eqnarray*}
Thus, $|q_{d+1}^{(1)}(z)| > |q_{d+1}^{(2)}(z)|$ holds for any
complex number $z$ with $\realpartof{z} \ge 0$.

Therefore, for any $d \ge 4$, the inequality
$|q_d^{(1)}(z)| > |q_d^{(2)}(z)|$ holds for any
complex number $z$ with a nonnegative real part.
This implies that the real part of any complex root of $i(\Pc_{K_d},m)$ is negative.

We shall also prove the other half, that $-\tfrac{d}{2} < \realpartof{\a}$.
To this end, 
it suffices to show that all roots of $j_d(l)=i\left(\Pc_{K_d}, -l-\frac{d}{2}\right)$ have negative real parts.
Let $r_d^{(1)}(l)$ and $r_d^{(2)}(l)$ be
\begin{eqnarray*}
  r_d^{(1)}(l) &=& {(-1)}^{d-1}q_d^{(1)}\left(-l-\frac{d}{2}\right)=(2l+1)\cdots(2l+d-1)\\
  r_d^{(2)}(l) &=& {(-1)}^{d-1}q_d^{(2)}\left(-l-\frac{d}{2}\right)=d\left(l-\frac{d-4}{2}\right)\cdots\left(l+\frac{d}{2}\right).
\end{eqnarray*}
Then for a complex number $z$, it holds that
\[
j_d(z)=0 \iff   r_d^{(1)}(z) =   r_d^{(2)}(z)
.\]
Let us prove $|r_d^{(1)}(z)| > |r_d^{(2)}(z)|$ for any
complex number $z$ with a nonnegative real part 
by mathematical induction on $d \geq 4$.

For $d=4$, it immediately follows from the inequality between $q_4^{(1)}$ and $q_4^{(2)}$:
\begin{eqnarray*}
|r_4^{(1)}(z)|\,=\,|q_4^{(1)}(z)| &>& |q_4^{(2)}(z)|\,=\,|r_4^{(2)}(z)| .
\end{eqnarray*}
And so we need $d=5$ also as a base case:
\begin{eqnarray*}
  |r_5^{(1)}(z)| &=& |2z+1||2z+2||2z+3||2z+4|\\
  &>& \frac{5}{4}  |z+1||2z+1||2z+3||2z+4|\\
  &>& \frac{5}{4}  |z-\tfrac{1}{2}||2z+1||2z+3||z+\tfrac{5}{2}|\\
  &=& 5 \left|z-\frac{1}{2}\right|\left|z+\frac{1}{2}\right|\left|z+\frac{3}{2}\right|\left|z+\frac{5}{2}\right|\\
  &=& |r_5^{(2)}(z)| .
\end{eqnarray*}

Assume for $d$ the validity of $|r_d^{(1)}(z)| > |r_d^{(2)}(z)|$ for any
complex number $z$ with $\realpartof{z} \geq 0$.

Then, from the fact that
\begin{eqnarray*}
  |r_{d+2}^{(1)}(z)|&=&|2z+d||2z+d+1| |r_{d}^{(1)}(z)|\\
  |r_{d+2}^{(2)}(z)|&=&\frac{d+2}{d}\left|z-\frac{d}{2}+1\right|\left|z+\frac{d}{2}+1\right| |r_{d}^{(2)}(z)|,
\end{eqnarray*}
it follows that
\begin{eqnarray*}
  d|r_{d+2}^{(1)}(z)| &=& d|2z+d||2z+d+1| |r_{d}^{(1)}(z)|\\
  &>& d |2z+d|\left|z+\tfrac{d}{2}+1\right| |r_{d}^{(2)}(z)|\\
  &=& |2dz+d^2|\left|z+\tfrac{d}{2}+1\right| |r_{d}^{(2)}(z)|\\
  &>& |(d+2)z+d^2-4|\left|z+\tfrac{d}{2}+1\right| |r_{d}^{(2)}(z)|\\
  &>& (d+2)\left|z-\frac{d-2}{2}\right|\left|z+\frac{d}{2}+1\right| |r_{d}^{(2)}(z)|\\
  &=& d|r_{d+2}^{(2)}(z)| .
\end{eqnarray*}
Thus, $|r_{d+2}^{(1)}(z)| > |r_{d+2}^{(2)}(z)|$ holds for any
complex number $z$ with $\realpartof{z} \ge 0$.

Therefore, for any $d \ge 4$, the inequality
$|r_d^{(1)}(z)| > |r_d^{(2)}(z)|$ holds for any
complex number $z$ with a nonnegative real part.
This implies that any complex root of $j_d(l)$ has a negative
real part.
\end{proof}

\subsection{Complete Multipartite Graphs}
\label{sec:epcmg}

We computed the roots of the Ehrhart polynomials $i(\Pc_G, m)$ of
complete multipartite graphs $G$ as well.
Since complete multipartite graphs are a special subclass of
connected simple graphs, our interest is mainly on the cases
where the general method could not complete the computation,
i.e., complete multipartite graphs of orders $d \ge 10$.

A complete multipartite graph of type $(q_1,\ldots,q_t)$, 
denoted by $K_{q_1,\ldots,q_t}$,
is constructed as follows.
Let $V(K_{q_1,\ldots,q_t}) = \Union_{i=1}^{t} V_i$ be a disjoint union of vertices
with $\numberof{V_i} = q_i$ for each $i$ 
and the edge set $E(K_{q_1,\ldots,q_t})$ 
be $\{\{u, v\}\ |\ u \in V_i,\ v \in V_j\ (i \neq j)\}$.
The graph $K_{q_1,\ldots,q_t}$ is unique up to isomorphism.

The Ehrhart polynomials for complete multipartite graphs are
explicitly given in~\cite{OH2000}:
\begin{equation}\label{ohformula}
i(\Pc_G, m) = \binom{d+2m-1}{d-1}-\sum_{k=1}^{t}\sum_{1\le i\le j \le q_k}\binom{j-i+m-1}{j-i}\binom{d-j+m-1}{d-j}
\end{equation}
where $d = \sum_{k=1}^{t} q_k$ is a partition of $d$ and $G=K_{q_1,\ldots,q_t}$.

Another simpler formula is newly obtained.
\begin{Proposition}\label{prop:fsum}
  The Ehrhart polynomial $i(\Pc_{G}, m)$
  of the edge polytope of a complete multipartite graph $G=K_{q_1,\ldots,q_t}$
  is 
  \[  i(\Pc_{G}, m) =   f(m; d, d) - \sum_{k=1}^{t} f(m; d, q_k),  \]
  where $d = \sum_{k=1}^{t} q_k$ and
  \[  f(m; d, j) = \sum_{k=1}^{j} p(m; d, k)  \]
  with
  \[  p(m; d, j) = \binom{j+m-1}{j-1} \binom{d-j+m-1}{d-j}.  \]
\end{Proposition}
\begin{proof}
  Let $G$ denote a complete multipartite graph $K_{q_1,\ldots,q_t}$.
  We start from the formula (\ref{ohformula}).

  First, it holds that
  \[ \binom{d+2m-1}{d-1}=f(m; d, d). \]
  On the one hand, $\binom{d+2m-1}{d-1}$ is the number of combinations
  with repetitions choosing $2m$ elements from a set of cardinality $d$.
  On the other hand, 
  \[ f(m;d,d) = \sum_{j=1}^{d} \binom{j+m-1}{j-1} \binom{d-j+m-1}{d-j} \]
  counts the same number of combinations as the sum of the number of
  combinations in which the $(m+1)$-th smallest number is $j$.

  Second, it holds that
  \[
  \sum_{k=1}^{t}\sum_{1\le i\le j \le q_k}\binom{j-i+m-1}{j-i}\binom{d-j+m-1}{d-j}
  = \sum_{k=1}^{t} f(m; d, q_k).
  \]
  Since the outermost summations are the same on both sides,
  it suffices to show that 
  \[
  \sum_{1\le i\le j \le q_k}\binom{j-i+m-1}{j-i}\binom{d-j+m-1}{d-j}
  = f(m; d, q_k).
  \]
  The summation of the left-hand side can be transformed as follows: 
  \begin{eqnarray*}
    &&\sum_{1\le i\le j \le q_k} \binom{j-i+m-1}{j-i}\binom{d-j+m-1}{d-j}\\
    &=& \sum_{j=1}^{q_k} \sum_{i=1}^{j} \binom{j-i+m-1}{j-i}\binom{d-j+m-1}{d-j}\\
    &=& \sum_{j=1}^{q_k} \binom{d-j+m-1}{d-j} \sum_{i=1}^{j} \binom{j-i+m-1}{j-i}\\
    &=& \sum_{j=1}^{q_k} \binom{d-j+m-1}{d-j} \binom{m+j-1}{j-1}\\
    &=& \sum_{j=1}^{q_k} p(m; d, j)\\
    &=& f(m; d, q_k)
  \end{eqnarray*}

  Finally, substituting these transformed terms into the original formula (\ref{ohformula}) 
  gives the desired result.
\end{proof}

\begin{figure}[htb!]
\centering%
\includegraphics[scale=0.4]{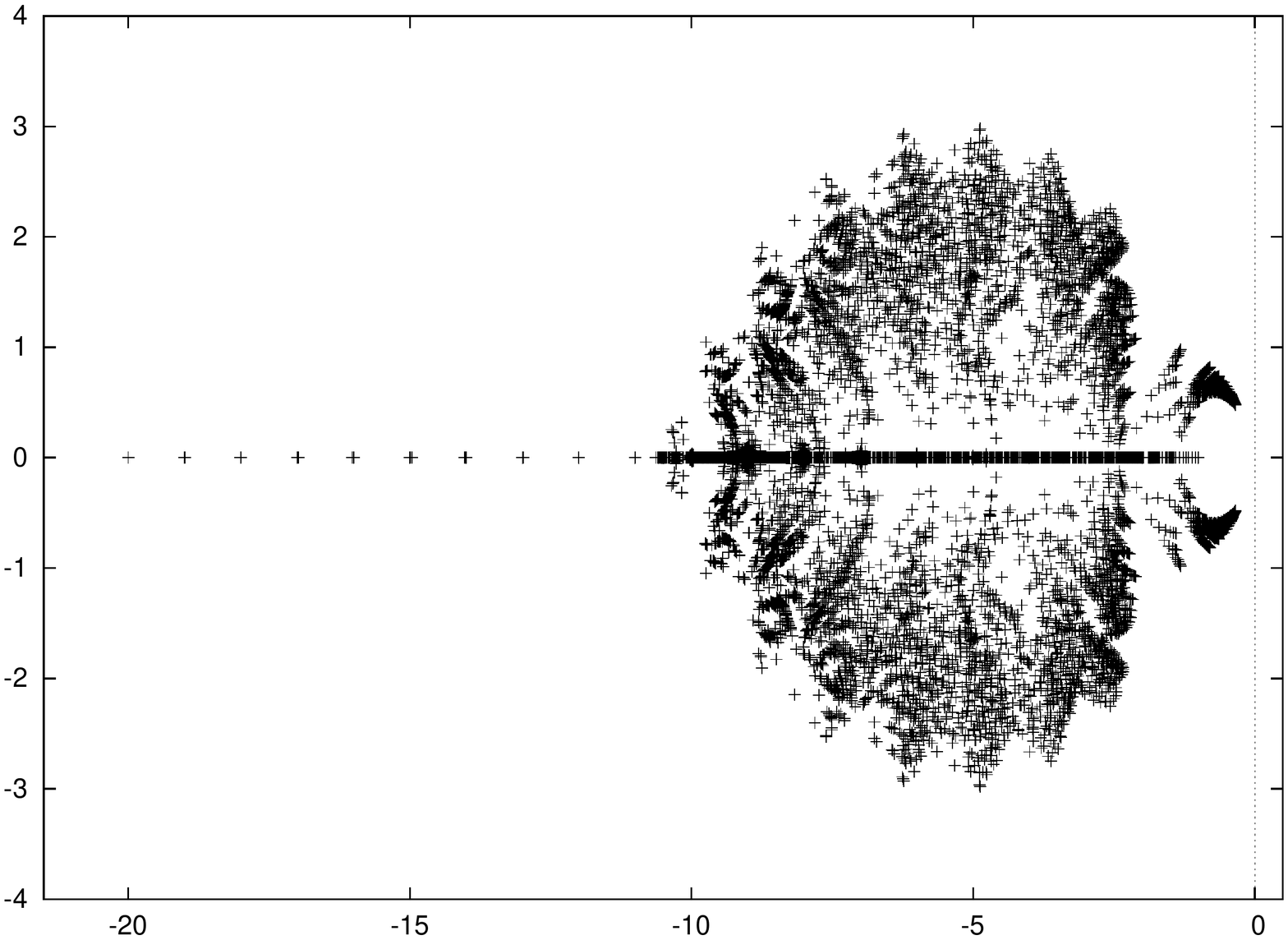}
\caption{$\Vb^{\text mp}_{22}$}
\label{fig:mp22}
\end{figure}

By the new formula above, we computed the roots of Ehrhart polynomials.
Let $\Vb^{\text mp}_d$ denote $\Union \Vb(i(\Pc_G, m))$, 
where the union runs over all complete multipartite graphs $G$ of order $d$.
Figure~\ref{fig:mp22} plots the points of $\Vb^{\text mp}_{22}$.
For all complete multipartite graphs of order 
$10$--$22$, Conjecture~\ref{conj:dstrip} holds.

Figure~\ref{fig:mp22}, for $\Vb^{\text mp}_{22}$, shows
that the noninteger roots lie in the circle
$\left|z + \tfrac{11}{2}\right| \le \tfrac{11}{2}$.
This fact is not exclusive to $22$ alone, but similar conditions 
hold for all $d \le 22$.
We conjecture:

\begin{Conjecture}\label{conj:circle}
  For any $d \ge 3$, 
  \[\Vb^{\text mp}_d  \subset
  \left\{z \in \CC\ |\ \left|z + \tfrac{d}{4}\right| \le \tfrac{d}{4}\right\}
  \union \{-(d-1),\ldots,-2,-1\}.\]
\end{Conjecture}

\begin{Remark}
(1) The leftmost point $-(d-1)$ can only be attained by $K_3$;
  this is shown in Proposition~\ref{prop:nodminusone}.
  Therefore, if we choose $d \ge 4$, the set of negative integers
  in the statement can be replaced with the set $\{-(d-2),\ldots,-2,-1\}$.
  However, $-(d-2)$ can be attained by the tree $K_{d-1,1}$
  for any $d$; see Example~\ref{eg:bipartite} below.

(2) Since \(0\) can never be a root of an Ehrhart polynomial,
  Conjecture~\ref{conj:circle} answers Question~\ref{q:leftplane}
  in the affirmative for complete multipartite graphs.
  Moreover, if Conjecture~\ref{conj:circle} holds, then
  Conjecture~\ref{conj:dstrip} holds for those graphs.

(3) The method of Pfeifle~\cite{Pfe2010} might be useful if the
  $\delta$-vector can be determined for edge polytopes of complete
  multipartite graphs.
\end{Remark}

\begin{Example}
\label{eg:bipartite}
The Ehrhart polynomial for complete bipartite graph $K_{p,q}$ is given in,
e.g.,~\cite[Corollary 2.7 (b)]{OH2000}:
\[ i(\Pc_{K_{p,q}}, m) = \binom{m+p-1}{p-1}\binom{m+q-1}{q-1},\]
and thus the roots are
\[ \Vb (i(\Pc_{K_{p,q}}, m)) = \{-1,\ldots, -\max(p-1, q-1)\} \]
and all of them are negative integers satisfying the condition in
Conjecture~\ref{conj:circle}.
\end{Example}

\begin{Example}
\label{eg:n11}
The edge polytope of a complete $3$-partite graph $\Pc_{K_{n,1,1}}$ for
$n \ge 2$ can be obtained as a pyramid from  $\Pc_{K_{n,2}}$ by adjoining
a vertex.
Therefore, its Ehrhart polynomial is the following:
\[ i(\Pc_{K_{n,1,1}}, m) = \sum_{j=0}^m i(\Pc_{K_{n,2}}, j). \]
Each term on the right-hand side is given in
Example~\ref{eg:bipartite} above.
By some elementary algebraic manipulations of binomial coefficients,
it becomes,
\[ i(\Pc_{K_{n,1,1}}, m) = \binom{m+n}{n}\frac{n m+n+1}{n+1}. \]
The noninteger root $\tfrac{-(n+1)}{n}$ is a real number in the circle
of Conjecture~\ref{conj:circle}.
\end{Example}

Now we prepare the following lemma 
for proving Proposition \ref{prop:nodminusone}. 

\begin{Lemma}\label{lem:psym}
  For any integer $1 \le j \le \tfrac{d}{2}$, the polynomial $p(m; d, j)$
  in Proposition~\ref{prop:fsum} satisfies:
  \[
  p(m; d, d-j) = \left(\frac{d}{j} - 1\right)p(m; d, j).
  \]
\end{Lemma}
\begin{proof}
  It is an easy transformation:
  \begin{eqnarray*}
    p(m; d, d-j) &=& \binom{(d-j)+m-1}{(d-j)-1} \binom{d-(d-j)+m-1}{d-(d-j)}\\
    &=& \binom{d-j+m-1}{d-j-1} \binom{j+m-1}{j}\\
    &=& \frac{d-j}{j} \binom{d-j+m-1}{d-j} \binom{j+m-1}{j-1}\\
    &=& \left(\frac{d}{j} - 1\right)p(m; d, j).
  \end{eqnarray*}
\end{proof}

\begin{Proposition}\label{prop:nodminusone}
  Let $(q_1,\ldots,q_t)$ be a partition of $d \ge 3$,
  satisfying $q_1 \ge q_2 \ge \cdots \ge q_t$.
  The Ehrhart polynomial $i(\Pc_{G}, m)$
  of the edge polytope of the complete multipartite graph $G=K_{q_1,\ldots,q_t}$
  does not have a root at $-(d-1)$
  except when the graph is $K_3$.
\end{Proposition}
\begin{proof}
  From Proposition~\ref{prop:fsum}, the Ehrhart polynomial of the edge polytope of
  $G=K_{q_1,\ldots,q_t}$ is
  \begin{eqnarray*}
    i(\Pc_{G}, m)
    &=& f(m; d, d) - \sum_{k=1}^{t} f(m; d, q_k)\\
    &=& p(m; d, d) + \sum_{j=1}^{d-1}p(m; d, j) - \sum_{k=1}^{t} \sum_{j=1}^{q_k} p(m; d, j)
  \end{eqnarray*}
  Since $p(m; d, d)$ has $-(d-1)$ as one of its roots,
  it suffices to show that the rest of the expression does not have $-(d-1)$
  as one of its roots.

  We evaluate $p(m; d, j)$ at $-(d-1)$ for $j$ from $1$ to $d-1$:
  \[
  p(-(d-1); d, j) = \binom{j-d}{j-1}\binom{-j}{d-j}
  \]
  by the definition of $p(m; d, j)$.
  If $j > 1$, its sign is ${(-1)}^{j-1 + d-j} = {(-1)}^{d-1}$
  since $j-d < 0$ and $-j < 0$.
  In case where $j=1$, since $j-1$ is zero, 
  \[
  p(-(d-1); d, 1) = \binom{-1}{d-1} = {(-1)}^{d-1}
  \]
  gives the same sign with other values of $j$.

  By the conjugate partition $(q'_1,\ldots,q'_{t'})$ of $(q_1,\ldots,q_t)$,
  which is given by $q'_j = \numberof{\{i \le t\ |\ q_i \ge j\}}$, we obtain
  \begin{eqnarray}\label{eq:conjpart}
    \sum_{j=1}^{d-1}p(m; d, j) - \sum_{k=1}^{t} \sum_{j=1}^{q_k} p(m; d, j)
    &=& \sum_{j=1}^{d-1}\left(1 - q'_j\right) p(m; d, j),
  \end{eqnarray}
  where we set, for simplicity, $q'_j = 0$ for $j > t'$.

  We show that all the coefficients of $p(m; d, j)$ are nonnegative
  for any $j$ from $1$ to $d - 1$ and there is at least one positive
  coefficient among them.\\
  {\bf (I)}  $q_1 \ge \tfrac{d}{2}$:\\
  The coefficients of $p(m; d, j)$ are zero for $q_1 \ge j \ge d-q_1$,
  unless $d = q_1 + q_2$, i.e., when the graph is a complete bipartite graph;
  the exceptional case will be discussed later.
  We assume, therefore, $q_2 < d - q_1$ for a while.
  Though equation (\ref{eq:conjpart}) gives the coefficient of $p(m; d, j)$
  as $1$ for $d>j>q_1$,
  by using Lemma~\ref{lem:psym}, we are able to let them be zero and
  the coefficient of $p(m; d, j)$ be $\tfrac{d}{j} - q'_j$ for $d-q_1>j>0$.
  Then all the coefficients of $p(m; d, j)$'s are positive,
  since the occurrence of integers greater than or equal to $j$
  in a partition of $d-q_1$ cannot be greater than $\tfrac{d-q_1}{j}$.\\
  {\bf (II)}  $q_1 < \tfrac{d}{2}$:\\
  Each coefficient of $p(m; d, j)$ in equation
  (\ref{eq:conjpart}) is $1$ for $d>j>\tfrac{d}{2}$.
  By Lemma~\ref{lem:psym}, we transfer them to lower $j$ terms
  so as to make the coefficients for $\tfrac{d}{2}>j>0$
  be $\tfrac{d}{j} - q'_j$.
  Then all the coefficients of $p(m; d, j)$'s are nonnegative,
  since the occurrence of integers greater than or equal to $j$
  in a partition of $d$ cannot be greater than $\tfrac{d}{j}$.
  Moreover, the coefficient is zero for at most one $j$,
  less than $\tfrac{d}{2}$.
  If $d=3$ and $q_1=q_2=q_3=1$, i.e., in case of $K_3$, there does not remain
  a positive coefficient. This exceptional case will be discussed later.

  For both (I) and (II), ignoring the exceptional cases, the terms on
  the right-hand side of equation (\ref{eq:conjpart})
  are all nonnegative when $d \equiv 1 \pmod 2$, or nonpositive otherwise,
  and there is at least one nonzero term.
  That is, $-(d-1)$ is not a root of 
  \[
  \sum_{j=1}^{d-1}p(m;d,j) - \sum_{k=1}^{t} \sum_{j=1}^{q_k} p(m;d,j).
  \]
  The Ehrhart polynomial $i(\Pc_{G}, m)$ is a sum of 
  a polynomial whose roots include $-(d-1)$ and
  another polynomial whose roots do not include $-(d-1)$.
  Therefore, $-(d-1)$ is not a root of $i(\Pc_{G}, m)$.

  Finally, we discuss the exceptional cases.
  The complete bipartite graphs are treated in Example~\ref{eg:bipartite}.
  In these cases, $-(d-1)$ is not a root of the Ehrhart polynomials.
  However, $-(d-1)=-2$ is actually a root of the Ehrhart polynomial
  of the edge polytope constructed from the complete graph $K_3$, as shown in
  Proposition~\ref{eg:complete} (1).
\end{proof}

\section{Edge polytopes of graphs with loops}
\label{sec:loops}


A convex polytope $\Pc$
of dimension $D$ is {\em simple}
if each vertex of $\Pc$
belongs to exactly $D$ edges of $\Pc$.
A simple polytope $\Pc$ is {\em smooth}
if at each vertex of $\Pc$,
the primitive edge directions form a lattice basis.

Now, if $e = \{ i, j \}$ is an edge of $G$,  
then $\rho(e)$ cannot be a vertex of $\Pc_G$ 
if and only if $i \neq j$ and $G$ has a loop
at each of the vertices $i$ and $j$. 
Suppose that $G$ has a loop at $i \in V(G)$ and $j \in V(G)$ and
that $\{i,j\}$ is not an edge of $G$.
Then $\Pc_G = \Pc_{G'}$ for the graph $G'$
defined by $E(G') = E(G) \cup \{ \{ i , j \} \}$.
Considering this fact, 
throughout this section, we assume that
$G$ satisfies the following condition:

\smallskip

\begin{enumerate}
\item[$(*)$]
If $i$, $j \in V(G)$ and 
if $G$ has a loop at each of $i$ and $j$, 
then the edge $\{ i, j \}$ belongs to $G$. 
\end{enumerate}

The graphs $G$ (allowing loops) whose edge polytope $\Pc_G$ is simple 
are completely classified by the following.
\begin{Theorem}\label{bunrui}{\em(\cite[Theorem 1.8]{OH2008})}
Let $W$ denote the set of vertices 
$i \in V(G)$
such that $G$ has no loop at $i$ and 
let $G'$ denote the induced subgraph of $G$ on $W$.
Then the following conditions are equivalent {\rm :}
\begin{itemize}
\item[(i)]
$\Pc_G$ is simple, but not a simplex {\rm ;}
\item[(ii)]
$\Pc_G$ is smooth, but not a simplex {\rm ;}
\item[(iii)]
$W \neq \emptyset$ and
$G$ is one of the following graphs {\rm :}
\begin{itemize}
\item[($\alpha $)]
$G$ is a complete bipartite graph with at least one cycle
of length $4$ {\rm ;}
\item[($\beta $)]
$G$ has exactly one loop, 
$G'$ is a complete bipartite graph and
if $G$ has a loop at $i$, then
$\{ i, j \} \in E(G)$ 
for all $j \in W$
{\rm ;}
\item[($\gamma $)]
$G$ has at least two loops,
$G'$ has no edge and
if $G$ has a loop at $i$, then
$\{ i, j \} \in E(G)$ 
for all $j \in W$.
\end{itemize}
\end{itemize}
\end{Theorem}


From the theory of Gr\"{o}bner bases, 
we obtain the Ehrhart polynomial $i(\Pc_G,m)$ of the edge polytope $\Pc_G$ above. 
In fact, 
\begin{Theorem}\label{normalizedvolume}{\em(\cite[Theorem 3.1]{OH2008})}
Let $G$ be a graph as in Theorem \ref{bunrui} (iii).
Let $W$ denote the set of vertices 
$i \in V(G)$
such that $G$ has no loop at $i$ and 
let $G'$ denote the induced subgraph of $G$ on $W$.
Then the Ehrhart polynomial $i(\Pc_G,m)$ of the edge polytope $\Pc_G$
are as follows:
\begin{itemize}
\item[($\alpha $)]
If $G$ is the complete bipartite graph 
on the vertex set 
$V_1 \cup V_2$ with $|V_1| = p$ and $|V_2| = q$,
then we have
\[
i(\Pc_G,m) = \binom{ p + m -1 }{ p -1 } \binom{ q + m -1 }{ q -1 }; 
\]
\item[($\beta $)]
If $G'$ is the complete bipartite graph 
on the vertex set 
$V_1 \cup V_2$ with $|V_1| = p$ and $|V_2| = q$,
then we have
\[
i(\Pc_G,m) = \binom{ p + m }{ p }\binom{ q + m }{ q }; 
\]
\item[($\gamma $)]
If $G$ possesses $p$ loops and $\numberof{V(G)} = d$, then we have
\[
i(\Pc_G,m) = 
\sum_{j=1}^p 
\binom{ j+m-2 }{ j-1}
\binom{ d-j+m }{ d-j}
.\]
\end{itemize}
\end{Theorem}

The goal of this section is to discuss the roots of Ehrhart polynomials of 
simple edge polytopes in Theorem \ref{bunrui} 
(Theorems \ref{root1}, \ref{root2}, and \ref{root3}).

\subsection{Roots of Ehrhart polynomials}

The consequences of the theorems above support Conjecture~\ref{conj:dstrip}.
Recall that ${\bf V}(f)$ denotes the set of roots of given polynomial $f$.


\begin{Example}
The Ehrhart polynomial for a graph $G$, 
the induced subgraph $G'$ of which is a complete bipartite graph $K_{p,q}$, 
is given in Theorem \ref{normalizedvolume} ($\beta$): 
\[
i(\Pc_G,n)=\binom{p+m}{p}\binom{q+m}{q}, 
\]
and thus the roots are 
\[
{\bf V}
\left(
\binom{ p + m }{ p } \binom{ q + m  }{ q }
\right)
=
\{
-1,-2,\ldots,-\max(p,q)
\}.
\]
\end{Example}

\begin{Example}
\label{gammarei}
Explicit computation of the roots of the Ehrhart polynomials 
obtained in Theorem \ref{normalizedvolume} ($\gamma$)
seems, in general, to be rather difficult. 

Let $p=2$.  Then
\begin{eqnarray*}
& & \binom{ m-1 }{ 0}
\binom{ d-1+m }{ d-1}
+
\binom{ m }{ 1}
\binom{ d-2+m }{ d-2}\\
&=&
\binom{ d-1+m }{ d-1}
+
m\binom{ d-2+m }{ d-2}\\
&=& 
\left(
\frac{d-1+m}{d-1}
+m
\right)
\binom{ d-2+m }{ d-2}
\\
&=& 
\frac{dm + d-1}{d-1}
\binom{ d-2+m }{ d-2}.
\end{eqnarray*}
Thus, 
\[
{\bf V}\left(i(\Pc_G,m)\right) = \left\{ -1,-2,\ldots,-(d-2), -\frac{d-1}{d} \right\}.
\]
Let $p=3$.  Then
\begin{eqnarray*}
& & 
\binom{ m-1 }{ 0}
\binom{ d-1+m }{ d-1}
+
\binom{ m }{ 1}
\binom{ d-2+m }{ d-2}+
\binom{ m +1 }{ 2}
\binom{ d-3+m }{ d-3}\\
&=&
\binom{ d-1+m }{ d-1}
+
m\binom{ d-2+m }{ d-2}
+\frac{m(m+1)}{2}
\binom{ d-3+m }{ d-3}\\
&=& 
\left(
\frac{(d-1+m)(d-2+m)}{(d-1)(d-2)}
+m
\frac{d-2+m}{d-2}
+\frac{m(m+1)}{2}
\right)
\binom{ d-3+m }{ d-3}
\end{eqnarray*}
and
\begin{eqnarray*}
& & 
\frac{(d-1+m)(d-2+m)}{(d-1)(d-2)}
+m
\frac{d-2+m}{d-2}
+\frac{m(m+1)}{2}
\\
& = &
\frac{2(d-1+m)(d-2+m)+2(d-1)m(d-2+m)+(d-1)(d-2)m(m+1)}{
2(d-1)(d-2)
}\\
& = &
\frac{(d^2-d+2) m^2+ (3d^2-5d) m+ (2d^2-6d+4)}{
2(d-1)(d-2)
}.
\end{eqnarray*}
Let
\[
f(m)=(d^2-d+2) m^2+ (3d^2-5d) m+ (2d^2-6d+4).
\]
Since $d > p =3$, one has
\begin{eqnarray*}
f(0) &= &2d^2-6d+4 =2(d-1)(d-2) >0;\\
f(-1) &=& (d^2-d+2) -(3d^2-5d) + (2d^2-6d+4)=-2d+6<0;\\
f(-2) &=& 4 (d^2-d+2) - 2(3d^2-5d) + (2d^2-6d+4)=12 > 0.
\end{eqnarray*}
Hence, 
\[
{\bf V}\left(i(\Pc_G,m)\right) = \{ -1,-2,\ldots,-(d-3), \alpha,\beta \}
\]
where $-2 < \alpha < -1 < \beta <0$.
\end{Example}

We try to find information about the roots of the 
Ehrhart polynomials obtained 
in Theorem \ref{normalizedvolume} ($\gamma$)  
with $d > p \geq 2$.

\begin{Theorem}
\label{root1}
Let $d$ and $p$ be integers with $d > p \geq 2$
and let
$$f_{d,p}(m) = \sum_{j=1}^p 
\binom{ j+m-2 }{ j-1}
\binom{ d-j+m }{ d-j}$$
be a polynomial of degree $d-1$ in the variable $m$.
Then 
$$
\{-1,-2,\ldots, -(d-p)\} \ \ \subset \ \ 
{\bf V} (f_{d,p}) \cap {\mathbb R} 
\ \ \subset \ \ 
[-(d-2),0).$$
\end{Theorem}

\begin{proof}
It is easy to see that
$f_{d,p}(0) = 1$ and
$f_{d,p}(m) > 0 $ for all $m>0 $.

From Example \ref{gammarei}, we may assume that $4 \leq p < d$.
Then 
\begin{eqnarray*}
& & f_{d,p}(m) \\
&=& 
\binom{ d-1+m }{ d-1}
+
m
\binom{ d-2+m }{ d-2}
+
\sum_{j=3}^p 
\binom{ j+m-2 }{ j-1}
\binom{ d-j+m }{ d-j}\\
&=& 
\left(
\frac{d-1+m}{d-1}+m
\right)
\binom{ d-2+m }{ d-2}
+
\sum_{j=3}^p 
\binom{ j+m-2 }{ j-1}
\binom{ d-j+m }{ d-j}\\
&=& 
\frac{md + d-1}{d-1}
\binom{ d-2+m }{ d-2}
+
\sum_{j=3}^p 
\binom{ j+m-2 }{ j-1}
\binom{ d-j+m }{ d-j}.
\end{eqnarray*}
If $m < -(d-2)$, then
$m+d-2 < 0$,  $md + d-1 < -(d-2) d + d-1 =  -(d-3)d -1  <0$,
\begin{eqnarray*}
m+d-j &\leq& m+d-3 <0\\
m+j-2 &\leq& m+p-2 \leq m+d-3 <0
\end{eqnarray*}
for each $j = 3, 4, \ldots,p$.
Hence, we have $(-1)^{d-1} f_{d,p}(m) >0$ for all $m < -(d-2)$.
Thus, we have
$
{\bf V} (f_{d,p}) \cap {\mathbb R} 
\ \subset \ 
[-(d-2),0)$.

Since
$$f_{d,p}(m) = 
\binom{ d-p+m }{ d-p}
\sum_{j=1}^p 
\binom{ j+m-2 }{ j-1}
\frac{
(d-j+m) \cdots (d-p+1+m)
}{
(d-j) \cdots (d-p+1)
}
,$$
it follows that
$$
{\bf V} \left(
\binom{ d-p+m }{ d-p}
\right) =
\{-1,-2,\ldots, -(d-p)\} \subset
{\bf V} (f_{d,p})
.
$$
\end{proof}

\begin{Theorem}
\label{root2}
Let $d$ and $p$ be integers with $d > p \geq 2$
and let $f_{d,p}(m)$ be the polynomial defined above. 
If $d -2p +2 \geq 0$, then 
$${\bf V} (f_{d,p}) = \{-1,-2,\ldots, -(d-p),\alpha_1,\alpha_2,\ldots, \alpha_{p-1}\}$$
where
\[
-(p-1) < \alpha_{p-1} < -(p-2) < \alpha_{p-2} < -(p-3) < \cdots < -1 < \alpha_1 < 0. 
\]
\end{Theorem}

\begin{proof}
Let
$$g_{d,p}(m) = 
\frac{
f_{d,p}(m)
}{\binom{ d-p+m }{ d-p}}
=
\sum_{j=1}^p 
\binom{ j+m-2 }{ j-1}
\frac{
(d-j+m) \cdots (d-p+1+m)
}{
(d-j) \cdots (d-p+1)
}.
$$
It is enough to show that 
\[
(-1)^k g_{d,p}(k) >0
\]
for $k = 0,-1,-2,\ldots, -(p-1)$.

\medskip

\noindent
{\bf (First Step)}
We claim that $(-1)^{-(p-1)} g_{d,p}(-(p-1)) >0$.
A routine computation on binomial coefficients yields 
the equalities
\begin{eqnarray*}
& &
g_{d,p}(-(p-1))\\
&=&
\frac{
\sum_{j=1}^p 
(-1)^{j-1}
\binom{p-1 }{ j-1}
\prod_{i=1}^{j-1} (d - i)
\prod_{k=j}^{p-1} (d - k - (p-1))
}{
(d-1) \cdots (d-p+1)
}
\end{eqnarray*}
and
\begin{eqnarray*}
& &
\sum_{j=1}^p 
(-1)^{j-1}
\binom{p-1 }{ j-1}
\prod_{i=1}^{j-1} (d - i)
\prod_{k=j}^{p-1} (d - k - (p-1))\\
&= & (-1)^{p-1} (p-1) p \cdots  (2 p- 3).
\end{eqnarray*}
Hence, 
\[
(-1)^{p-1}  g_{d,p}(-(p-1)) 
= \frac{(p-1) p \cdots  (2 p- 3)}{(d-1) \cdots (d-p+1)} > 0.
\]

\medskip

\noindent
{\bf (Second Step)}
Working with induction on $p$,
we now show that
\[
(-1)^k g_{d,p}(k) >0
\]
for $k = 0,-1,-2,\ldots, -(p-2)$.
Again, a routine computation on binomial coefficients yields
\begin{eqnarray*}
g_{d,p}(m)
=
\binom{ p+m-2 }{ p-1}
+
\frac{d-p+1+m}{d-p+1}
g_{d,p-1}(m).
\end{eqnarray*}
Hence, 
$$(-1)^k g_{d,p}(k) =  \frac{d-p+1+k}{d-p+1}
(-1)^k g_{d,p-1}(k).$$
Since $d - 2p +2 \geq 0$, one has
$$d-p+1+k \geq d-p+1-(p-2) = d-2p +3 > 0.$$
By virtue of $(-1)^{-(p-1)} g_{d,p}(-(p-1)) >0$, 
together with
the hypothesis of induction, 
it follows that \[
(-1)^k g_{d,p-1}(k) >0.
\]
Thus, 
\[
(-1)^k g_{d,p}(k) >0,
\] 
as desired.
\end{proof}

If $d -2p +2 \geq 0$, then
it follows that
$$
\floorof{ \frac{d-1}{2} } \leq d-p
.$$
In this case, around half of the elements of ${\bf V} (f_{d,p})$ 
are negative integers.
This fact remains true even if $d -2p +2 < 0$.

\begin{Theorem}
\label{root3}
Let $d$ and $p$ be integers with $d > p \geq 2$
and let $f_{d,p}(m)$ be the polynomial defined above. 
Then 
$$
\left\{-1,-2,\ldots, - \left\lfloor \frac{d-1}{2} \right\rfloor \right\}
\subset
{\bf V} (f_{d,p}).$$
\end{Theorem}

\begin{proof}
If $d -2p +2 \geq 0$, then 
it follows from Theorem \ref{root1}.
(Note that if $p=2$, then $d -2p +2 = d-2 >0$.)

Work with induction on $p$.
Let $d -2p +2 < 0$.
By Theorem \ref{root1},
it is enough to show that 
$g_{d,p}(k) = 0$ for all $k=-(d-p+1),\ldots, -\left\lfloor \frac{d-1}{2} \right\rfloor$.
As in the proof of Theorem \ref{root2},
we have
$$g_{d,p}(m)
=
\binom{ p+m-2 }{ p-1}
+
\frac{d-p+1+m}{d-p+1}
g_{d,p-1}(m).
$$
Since $d -2p +2 < 0$, it follows that $ \left\lfloor \frac{d-1}{2} \right\rfloor \leq p-2$.
Thus, 
$$g_{d,p}(k)
=
\frac{d-p+1+k}{d-p+1}
g_{d,p-1}(k).
$$
By virtue of
$$g_{d,p}(-(d-p+1))
=
\frac{0}{d-p+1}
g_{d,p-1}(-(d-p+1))
=0
$$
together with the hypothesis of induction, 
it follows that
$g_{d,p}(k) = 0$ for all $k=-(d-p+1),\ldots, -\left\lfloor \frac{d-1}{2} \right\rfloor$.
\end{proof}

\begin{Example}
Let 
$d =12$.
Then $d -2p +2 \geq 0$ if and only if $p \leq 7$.
For $p=2,3,\ldots,7$, the roots of the Ehrhart polynomials are
$-1,-2,\ldots, -(d-p)=p -12$, together with the real numbers 
listed as follows:
$$
\begin{array}{lcccccc}
p=2  & -0.92 &  & & & &  \\
p=3  & -1.92 & -0.85 & & & &  \\
p=4  & -2.90 & -1.83 & -0.80  & & &  \\
p=5  & -3.83 & -2.77 & -1.74 & -0.76 & &  \\
p=6  & -4.67 & -3.65 & -2.65 & -1.66 & -0.72 &  \\
p=7  & -5.31 & -4.42 & -3.47 & -2.53 & -1.58 & -0.69  \\
\end{array}
$$
For $p=8,9,10,11$, the roots of the Ehrhart polynomials are
$-1,-2,-3,-4,-5= -\left\lfloor \frac{d-1}{2} \right\rfloor$, 
together with the following complex numbers:
$$
\small
\begin{array}{lcccccccc}
p=8  & -5.56 &  -4.19 & -3.31 & -2.41 & -1.51 & -0.65  \\
p=9  & -5.47 &  -4.79 & -3.16 & -2.29 & -1.43 & -0.62  \\
p=10  & -5.51 & -4.16 + 0.18 i &-4.16 - 0.18 i  & -2.16 & -1.34 & -0.59   \\
p=11 & -5.50 & -4.53 &  -3.08 + 0.06 i & -3.08 - 0.06 i & -1.24 & -0.55 
\end{array}
$$
(Computed by {\tt Maxima}~\cite{maxima}) 
Thus, in particular, the real parts of all roots are negative.
\end{Example}

\section{Symmetric Edge Polytopes}
\label{sec:symmetric}


Among the many topics explored in recent papers 
on the roots of Ehrhart polynomials of convex polytopes, 
one of the most fascinating is the Gorenstein Fano polytope.

Let $\Pc \subset \RR^d$ be an integral convex polytope 
of dimension $d$.
\begin{itemize}
\item
We say that $\Pc$ is a {\em Fano polytope} if 
the origin of $\RR^d$ is the unique integer point 
belonging to the interior of $\Pc$. 
\item
A Fano polytope is said to be {\em Gorenstein} 
if its dual polytope is integral.
(Recall that the dual polytope $\Pc^\vee$
of a Fano polytope $\Pc$ is a convex polytope
that consists of those $x \in \RR^d$
such that $\langle x, y \rangle \leq 1$ for all
$y \in \Pc$, where $\langle x, y \rangle$
is the usual inner product of $\RR^d$.)
\end{itemize}

In this section, we will prove that 
symmetric edge polytopes arising from finite connected simple graphs 
are Gorenstein Fano polytopes (Proposition \ref{GFp}). 
Moreover, we will consider the condition of unimodular equivalence 
(Theorem \ref{nonequivalent}). 
In addition, we will compute the Ehrhart polynomials of symmetric edge polytopes 
and discuss their roots. 

\subsection{Fano polytopes arising from graphs}

Throughout this section, 
let $G$ denote a simple graph on the vertex set 
$V(G) = \{ 1, \ldots, d \}$ with $E(G) = \{e_1,\ldots,e_n\}$ being the edge set.
Moreover, let $\Pc_G^{\pm} \subset \RR^d$ denote a symmetric
edge polytope constructed from $G$.


Let $\Hc \subset \RR^d$ denote the hyperplane defined by the equation 
$x_1+x_2+\cdots+x_d=0$. 
Now, since the integral points $\pm\sigma(e_1),\ldots,\pm\sigma(e_n)$ 
lie on the hyperplane $\Hc$, 
we have $\dim(\Pc_G^{\pm}) \leq d-1$. 

\begin{Proposition}
One has $\dim(\Pc_G^{\pm})=d-1$ if and only if 
$G$ is connected. 
\end{Proposition}
\begin{proof}
Suppose that $G$ is not connected. 
Let $G_1,\ldots,G_m$ with $m>1$ denote the connected components 
of $G$. Let, say, $\{1,\ldots,d_1\}$ be the vertex set of $G_1$ 
and $\{d_1+1,\ldots,d_2\}$ the vertex set of $G_2$. 
Then $\Pc_G^{\pm}$ lies on two hyperplanes defined by the equations 
$x_1+\cdots+x_{d_1}=0$ and $x_{d_1+1}+\cdots+x_{d_2}=0$. 
Thus, $\dim(\Pc_G^{\pm}) < d-1.$ 

Next, we assume that $G$ is connected. 
Suppose that $\Pc_G^{\pm}$ lies on the hyperplane defined by the equation 
$a_1x_1+\cdots+a_dx_d=b$ with $a_1,\ldots,a_d,b \in \ZZ$. 
Let $e=\{i,j\}$ be an edge of $G$. 
Then because $\sigma(e)$ lies on this hyperplane together with $-\sigma(e)$, 
we obtain $$a_i-a_j=-(a_i-a_j)=b.$$ Thus $a_i=a_j$ and $b=0$. 
For all edges of $G$, since $G$ is connected, 
we have $a_1=a_2=\cdots=a_d$ and $b=0$. 
Therefore, $\Pc_G^{\pm}$ lies only on the hyperplane $x_1+x_2+\cdots+x_d=0$. 
\end{proof}

For the rest of this section, we assume that $G$ is connected. 

\begin{Proposition}\label{GFp}
Let $\Pc_G^{\pm}$ be a symmetric edge polytope of a finite graph $G$. 
Then $\Pc_G^{\pm} \subset \Hc$ is a Gorenstein Fano polytope of dimension $d-1$. 
\end{Proposition}
\begin{proof}
Let $\phi:\RR^{d-1} \rightarrow \Hc$ be the bijective homomorphism with 
$$\phi(y_1,\ldots,y_{d-1})=(y_1,\ldots,y_{d-1},-(y_1+\cdots+y_{d-1})).$$ 
Thus, we can identify $\Hc$ with $\RR^{d-1}$. 
Therefore, $\phi^{-1}(\Pc_G^{\pm})$ is isomorphic to $\Pc_G^{\pm}$. 

Since one has 
$$
\frac{1}{2n}\sum_{j=1}^n\sigma(e_j) + \frac{1}{2n}\sum_{j=1}^n(-\sigma(e_j)) 
= (0,\ldots,0) \in \RR^d, 
$$
the origin of $\RR^d$ is contained in 
the relative interior of $\Pc_G^{\pm} \subset \Hc$. 
Moreover, since 
$$
\Pc_G^{\pm} \subset \{(x_1,\ldots,x_d) \in \RR^d\ |\ -1 \leq x_i \leq 1,i=1,\ldots,d \}, 
$$
it is not possible for an integral point 
to exist anywhere in the interior of $\Pc_G^{\pm}$ except at the origin. 
Thus, $\Pc_G^{\pm} \subset \Hc$ is a Fano polytope of dimension $d-1$. 

%

Next, we prove that $\Pc_G^{\pm}$ is Gorenstein. 
Let $M$ be an integer matrix whose row vectors are 
$\sigma(e)$ or $-\sigma(e)$ with $e \in E(G)$. 
Then $M$ is a totally unimodular matrix. 
From the theory of totally unimodular matrices (\cite[Chapter 9]{Schrijver1986}), 
it follows that a system of equations $yA=(1,\ldots,1)$ 
has integral solutions, where $A$ is a submatrix of $M$. 
This implies that the equation of each supporting 
hyperplane of $\Pc_G^{\pm}$ is of the form 
$a_1x_1 + \cdots + a_dx_d = 1$ with each $a_i \in \ZZ$. 
In other words, the dual polytope of $\Pc_G^{\pm}$ is integral. 
Hence, $\Pc_G^{\pm}$ is Gorenstein, as required. 
\end{proof}

\subsection{When is $\Pc_G^{\pm}$ unimodular equivalent?}

In this subsection, we consider the conditions under which 
$\Pc_G^{\pm}$ is unimodular equivalent with 
$\Pc_{G'}^{\pm}$ for graphs $G$ and $G'$. 

Recall that for a connected graph $G$, we call $G$ a {\em 2-connected} graph
if the induced subgraph with the vertex set $V(G) \backslash \{i\}$ 
is still connected for any vertex $i$ of $G$.

Let us say a Fano polytope $\Pc \subset \RR^d$ {\em splits} into 
$\Pc_1$ and $\Pc_2$ if $\Pc$ is the convex hull of the two Fano polytopes 
$\Pc_1 \subset \RR^{d_1}$ and $\Pc_2 \subset \RR^{d_2}$ with $d=d_1+d_2$. 
That is, by arranging the numbering of coordinates, we have 
$$\Pc = \con(\{ (\a_1,{\bf 0}) \in \RR^d\ |\ \a_1 \in \Pc_1\} 
\union \{({\bf 0},\a_2) \in \RR^d\ |\ \a_2 \in \Pc_2 \}).$$ 

\begin{Lemma}\label{split}
$\Pc_G^{\pm}$ cannot split if and only if $G$ is 2-connected. 
\end{Lemma}
\begin{proof}
{\bf (``Only if'')} Suppose that $G$ is not 2-connected, 
i.e., there is a vertex $i$ of $G$ such that 
the induced subgraph $G'$ of $G$ with the vertex set 
$V(G) \backslash \{i\}$ is not connected. 
For a matrix 
\begin{eqnarray}\label{vmat}
\begin{pmatrix}
\sigma(e_1)  \\
-\sigma(e_1) \\
\vdots     \\
\sigma(e_n)  \\
-\sigma(e_n) \\
\end{pmatrix}
\end{eqnarray}
whose row vectors are the vertices of $\Pc_G^{\pm}$, 
we add all the columns of (\ref{vmat}) except the $i$\th column 
to the $i$\th column. 
Then the $i$\th column vector becomes equal to the zero vector. 
Let, say, $\{1,\ldots,i-1\}$ and $\{i+1,\ldots,d\}$ 
denote the vertex set of the connected components of $G'$. 
Then, by arranging the row vectors of (\ref{vmat}) if necessary, 
the matrix (\ref{vmat}) can be transformed into 
\begin{eqnarray*}
\begin{pmatrix}
M_1 &0 \\
0   &M_2
\end{pmatrix}. 
\end{eqnarray*}
This means that $\Pc_G^{\pm}$ splits into $\Pc_1$ and $\Pc_2$, 
where the vertex set of $\Pc_1$ (respectively $\Pc_2$) constitutes 
the row vectors of $M_1$ (respectively $M_2$). 

{\bf (``If'')} We assume that $G$ is 2-connected. 
Suppose that $\Pc_G^{\pm}$ splits into $\Pc_1,\ldots,\Pc_m$ 
and each $\Pc_i$ cannot split, where $m>1$. 
Then by arranging the row vectors if necessary, 
the matrix (\ref{vmat}) can be transformed into 
\begin{eqnarray*}
\begin{pmatrix}
&M_1            &         &\text{\Large 0} \\
&               &\ddots   &           \\
&\text{\Large 0}&         &M_m \\
\end{pmatrix}. 
\end{eqnarray*}
Now, for a row vector $v$ of each matrix $M_i$, 
$-v$ is also a row vector of $M_i$. Let 
$$
v_{i_1},\ldots,v_{i_{k_i}},-v_{i_1},\ldots,-v_{i_{k_i}}
$$
denote the row vectors of $M_i$, where $e_{i_1},\ldots,e_{i_{k_i}}$ are 
the edges of $G$ with $v_{i_j}=\sigma(e_{i_j})$ or $v_{i_j}=-\sigma(e_{i_j})$, 
and $G_i$ denote the subgraph of $G$ with the edge set $\{e_{i_1},\ldots,e_{i_{k_i}}\}$. 
Then for the subgraphs $G_1,\ldots,G_m$ of $G$, one has 
\begin{eqnarray}\label{ineq}
|V(G_1)| + \cdots + |V(G_m)| \geq d+2(m-1), 
\end{eqnarray}
where $V(G_i)$ is the vertex set of $G_i$. \\
(In fact, the inequality (\ref{ineq}) follows by induction on $m$. 
When $m=2$, since $G$ is 2-connected, 
$G_1$ and $G_2$ share at least two vertices. 
Thus, one has $|V(G_1)|+|V(G_2)| \geq d+2$. 
When $m=k+1$, since $G$ is 2-connected, one has 
$$
\left| (\union_{i=1}^kV(G_i)) \cap V(G_{k+1}) \right| \geq 2. 
$$
Let $d'$ be the sum of the numbers of the columns of 
$M_1,\ldots,M_{k-1}$ and $M_k$ 
and $d''$ be the number of the columns of $M_{k+1}$, 
where $d'+d''=d$. Then one has 
\begin{eqnarray*}
|V(G_1)| + \cdots + |V(G_k)| + |V(G_{k+1})| &\geq& 
d'+ 2(k-1) +|V(G_{k+1})| \\
&\geq& d'+d''+2(k-1)+2 = d+2k 
\end{eqnarray*} 
by the hypothesis of induction.) \\
In addition, each $\Pc_{G_i}^{\pm}$ cannot split. 
Thus one has $\dim(\Pc_{G_i}^{\pm})=|V(G_i)|-1$ 
since each $G_i$ is connected by the proof of the ``only if'' part. 
It then follows from this equality and the inequality (\ref{ineq}) that 
\begin{eqnarray*}
d-1 &=& \dim(\Pc_{G_1}^{\pm}) + \cdots + \dim(\Pc_{G_m}^{\pm}) 
= |V(G_1)| + \cdots + |V(G_m)| - m \\
&\geq& d+2m-2-m =d+m-2 \geq d \quad\quad (m \geq 2), 
\end{eqnarray*}
a contradiction. Therefore, $\Pc_G^{\pm}$ cannot split. 
\end{proof}

\begin{Lemma}\label{twoconnected}
Let $G$ be a 2-connected graph. Then, for a graph $G'$, 
$\Pc_G^{\pm}$ is unimodular equivalent with $\Pc_{G'}^{\pm}$ 
as an integral convex polytope if and only if 
$G$ is isomorphic to $G'$ as a graph. 
\end{Lemma}
\begin{proof}
If $|V(G)|=2$, the statement is obvious. Thus, we assume that $|V(G)|>2$. \\
{\bf (``Only if'')} 
Suppose that $\Pc_G^{\pm}$ is unimodular equivalent with $\Pc_{G'}^{\pm}$. 
Let $M_G$ (respectively $M_{G'}$) denote the matrix whose row vectors are the vertices of 
$\Pc_G^{\pm}$ (respectively $\Pc_{G'}^{\pm}$). 
Then there is a unimodular transformation $U$ such that one has 
\begin{eqnarray}\label{unimodular}
M_GU=M_{G'}. 
\end{eqnarray}
Thus, each row vector of $M_G$, i.e., each edge of $G$, 
one-to-one corresponds to each edge of $G'$. 
Hence, $G$ and $G'$ have the same number of edges. 
Moreover, since $G$ is 2-connected, 
$\Pc_G^{\pm}$ cannot split by Lemma \ref{split}. 
Thus, $\Pc_{G'}^{\pm}$ also cannot split; that is to say, 
$G'$ is also 2-connected. In addition, 
if we suppose that $G$ and $G'$ do not have the same number of vertices, 
then $\dim(\Pc_G^{\pm}) \not= \dim(\Pc_{G'}^{\pm})$ since $G$ and $G'$ are connected, 
a contradiction. Thus, the number of the vertices of $G$ 
is equal to that of $G'$. 

Now an arbitrary 2-connected graph with $|V(G)|>2$ 
can be obtained by the following method: 
start from a cycle and 
repeatedly append an $H$-path to a graph $H$ that has been already constructed. 
(Consult, e.g., \cite{Wilson}.) 
In other words, there is one cycle $C_1$ and $(m-1)$ paths 
$\Gamma_2,\ldots,\Gamma_m$ such that 
\begin{eqnarray}\label{const}
G=C_1 \union \Gamma_2 \union \cdots \union \Gamma_m. 
\end{eqnarray}
Under the assumption that 
$G$ is 2-connected and one has the equality (\ref{unimodular}), 
we show that $G$ is isomorphic to $G'$ by induction on $m$. 

If $m=1$, i.e., $G$ is a cycle, then $G$ has $d$ edges. 
Let $a_i,i=1,\ldots,d$ denote the degree of each vertex $i$ of $G'$. 
Then one has 
$$
a_1+a_2+\cdots+a_d=2d. 
$$
If there is $i$ with $a_i=1$, then $G'$ is not 2-connected. 
Thus, $a_i \geq 2$ for $i=1,\ldots,d$. Hence, $a_1=\cdots=a_d=2.$ 
It then follows that $G'$ is also a cycle of the same length as $G$, 
which implies that $G$ is isomorphic to $G'$. 

When $m=k+1$, we assume (\ref{const}). Let $\tilde{G}$ denote 
the subgraph of $G$ with 
$$\tilde{G} = C_1 \union \Gamma_2 \union \cdots \union \Gamma_k.$$ 
Then $\tilde{G}$ is a 2-connected graph. 
Since each edge of $G$ has one-to-one correspondence with each edge of $G'$, 
there is a subgraph $\tilde{G'}$ of $G'$ 
each of whose edges corresponds to those of $\tilde{G}$. 
Then one has $M_{\tilde{G}}U=M_{\tilde{G'}}$, 
where $M_{\tilde{G}}$ (respectively $M_{\tilde{G'}}$) is a submatrix of 
$M_G$ (respectively $M_{G'}$) whose row vectors are the vertices of 
$\Pc_{\tilde{G}}^{\pm}$ (respectively $\Pc_{\tilde{G'}}^{\pm}$). 
Thus, $\tilde{G}$ is isomorphic to $\tilde{G'}$ 
by the hypothesis of induction. Let 
$\Gamma_{k+1}=(i_0,i_1,\ldots,i_p)$ with $i_0 < i_1 < \cdots < i_p$ 
and $e_{i_l}=\{i_{l-1},i_l\},l=1,\ldots,p$ denote the edges of $\Gamma_{k+1}$. 
In addition, let $e_{i_1}',\ldots,e_{i_p}'$ denote the edges of $G'$ 
corresponding to the edges $e_{i_1},\ldots,e_{i_p}$ of $G$. 
Here, the edges $e_{i_1}',\ldots,e_{i_p}'$ of $G'$ 
are not the edges of $\tilde{G'}$. 
Since $i_0$ and $i_p$ are distinct vertices of $\tilde{G}$ 
and $\tilde{G}$ is connected, there is a path 
$\Gamma=(i_0,j_1,j_2,\ldots,j_{q-1},i_p)$ with 
$i_0=j_0 < j_1 < j_2 < \cdots < j_{q-1} < j_q=i_p$ in $\tilde{G}$. 
Let $e_{j_l}=\{j_{l-1},j_l\},l=1,\ldots,q$ denote the edges of $\Gamma$. 
Then by renumbering the vertices of $\tilde{G'}$ if necessary, 
there is a path $\Gamma'=(i_0',j_1',j_2',\ldots,j_{q-1}',i_p')$ 
with $i_0'=j_0' < j_1' < j_2' < \cdots < j_{q-1}' < j_q'=i_p'$ in $\tilde{G'}$ 
since $\tilde{G}$ is isomorphic to $\tilde{G'}$. 
Let $e_{j_l}'=\{j_{l-1}',j_l'\},l=1,\ldots,q$ denote the edges of $\Gamma'$. 
However, by (\ref{unimodular}), 
each edge $e_{j_l}$ of $\tilde{G}$ has one-to-one correspondence with 
each edge $e_{j_l}''$ of $\tilde{G'}$. 
Thus, each edge $e_{j_l}'$ of $\tilde{G'}$ has one-to-one correspondence with 
each edge $e_{j_l}''$ of $\tilde{G'}$. In other words, one has 
$$
\{e_{j_l}'\ |\ l=1,\ldots,q\} = \{e_{j_l}''\ |\ l=1,\ldots,q\}. 
$$

Since there are $\Gamma_{k+1}$ and $\Gamma$ that are paths 
from $i_0$ to $i_p$, one has 
\begin{eqnarray}\label{equation}
\sum_{l=1}^p\sigma(e_{i_l}) = \sum_{l=1}^q\sigma(e_{j_l}). 
\end{eqnarray}
On the one hand, if we multiply the left-hand side of the equation (\ref{equation}) with $U$, 
then we have 
$$
\sum_{l=1}^p\sigma(e_{i_l})U = \sum_{l=1}^p\sigma(e_{i_l}'). 
$$
On the other hand, if we multiply the right-hand side of the equation (\ref{equation}) with $U$, 
then we have 
$$
\sum_{l=1}^q\sigma(e_{j_l})U = \sum_{l=1}^q\sigma(e_{j_l}'') = 
\sum_{l=1}^q\sigma(e_{j_l}') = {\bf e}_{i_0'} - {\bf e}_{i_p'}. 
$$
Hence, we have $\sum_{l=1}^p\sigma(e_{i_l}')={\bf e}_{i_0'} - {\bf e}_{i_p'}$. 
This means that the edges $e_{i_1}',\ldots,e_{i_p}'$ of $G'$ 
construct a path from the vertex $i_0'$ to $i_p'$, 
which is isomorphic to $\Gamma_{k+1}$. 
Therefore, $G$ is isomorphic to $G'$. \\
{\bf (`` if '')} 
Suppose that $G$ is isomorphic to $G'$. 
Then by renumbering the vertices if necessary, 
it can be easily verified that $\Pc_G^{\pm}$ is 
unimodular equivalent with $\Pc_{G'}^{\pm}$. 
\end{proof}

\begin{Theorem}\label{nonequivalent}
For a connected simple graph $G$ (respectively $G'$), 
let $G_1,\ldots,G_m$ (respectively $G_1',\ldots,G_{m'}'$) 
denote the 2-connected components of $G$ (respectively $G'$). 
Then $\Pc_G^{\pm}$ is unimodular equivalent with $\Pc_{G'}^{\pm}$ 
if and only if $m=m'$ and $G_i$ is isomorphic to $G_i'$ 
by renumbering if necessary. 
\end{Theorem}
\begin{proof}
It is clear from Lemma \ref{split} and Lemma \ref{twoconnected}. 
If $G_i$ is isomorphic to $G_i'$ for $i=1,\ldots,m$, 
by virtue of Lemma \ref{split} and Lemma \ref{twoconnected}, 
then $\Pc_G^{\pm}$ is unimodular equivalent with $\Pc_{G'}^{\pm}$. 
On the contrary, suppose that $\Pc_G^{\pm}$ is unimodular equivalent with $\Pc_{G'}^{\pm}$. 
If $m \not= m'$, one has a contradiction by Lemma \ref{split}. Thus, $m=m'$. 
Moreover, by our assumption, $G_i$ is isomorphic to $G_i'$ by Lemma \ref{twoconnected}. 
\end{proof}

\subsection{Roots of the Ehrhart polynomials of $\Pc_G^{\pm}$} 
In this subsection, we study the Ehrhart polynomials of $\Pc_G^{\pm}$ 
and their roots. 

Let $\Pc \subset \RR^d$ be a Fano polytope with 
$\delta(\Pc) = (\delta_0, \delta_1, \ldots, \delta_d)$ being 
its $\delta$-vector. 
It follows from \cite{Bat1994} and \cite{Hib1992} 
that the following conditions are equivalent: 
\begin{itemize}
\item
$\Pc$ is Gorenstein; 
\item
$\delta(\Pc)$ is symmetric, i.e., 
$\delta_i = \delta_{d-i}$ for every $0 \leq i \leq d$; 
\item
$i(\Pc, m) = ( - 1 )^{d} i(\Pc, - m - 1)$. 
\end{itemize}
Since $i(\Pc, m) = ( - 1 )^{d} i(\Pc, - m - 1)$, 
the roots of $i(\Pc, m)$ locate symmetrically 
in the complex plane with respect to the line $\Re(z) = - \tfrac{1}{2}$.

\begin{Proposition}
If $G$ is a tree, then $\Pc_G^{\pm}$ is unimodular equivalent with 
\begin{eqnarray}\label{crosspolytope}
\con(\{\pm{\bf e}_1,\ldots,\pm{\bf e}_{d-1}\}). 
\end{eqnarray}
\end{Proposition}
\begin{proof}
If $G$ is a tree, 
then any 2-connected component of $G$ consists of one edge 
and $G$ possesses $(d-1)$ 2-connected components. 
Thus, by Theorem \ref{nonequivalent}, for any tree $G$, 
$\Pc_G^{\pm}$ is unimodular equivalent. 
Hence we should prove only the case where $G$ is a path, i.e., 
the edge set of $G$ is $\{\{i,i+1\}\ |\ i=1,\ldots,d-1\}$. 

Let 
\begin{eqnarray*}
\begin{pmatrix}
\sigma(e_1)  \\
-\sigma(e_1) \\
\vdots     \\
\sigma(e_{d-1})  \\
-\sigma(e_{d-1}) \\
\end{pmatrix}
\end{eqnarray*}
denote the matrix whose row vectors are the vertices of $\Pc_G^{\pm}$, 
where $e_i=\{i,i+1\},i=1,\ldots,d-1$ are the edges of $G$. 
If we add the $d$\th column to the $(d-1)$\th column, 
the $(d-1)$\th column to the $(d-2)$\th column, $\ldots$, and 
the second column to the first column, 
then the above matrix is transformed into 
\begin{eqnarray*}
\begin{pmatrix}
&0      &M        &         &{\bf 0} \\
&\vdots &         &\ddots   &        \\
&0      &{\bf 0}  &         &M       \\
\end{pmatrix}, 
\end{eqnarray*}
where $M$ is the $2 \times 1$ matrix 
$\begin{pmatrix}
-1 \\
1
\end{pmatrix}$. This implies that $\Pc_G^{\pm}$ is 
unimodular equivalent with (\ref{crosspolytope}). 
\end{proof}

Let $(\delta_0,\delta_1,\ldots,\delta_{d-1}) \in \ZZ^d$ be 
the $\delta$-vector of (\ref{crosspolytope}). 
Then it can be calculated that 
$$
\delta_i=\binom{d-1}{i}, i=0,1,\ldots,d-1. 
$$
It then follows from the well-known theorem \cite{Rod2002} that if $G$ is tree, 
the real parts of all the roots of $i(\Pc_G^{\pm}, m)$ are equal to $-\tfrac{1}{2}$. 
That is to say, all the roots $z$ of $i(\Pc_G^{\pm}, m)$ lie on the 
vertical line $\Re(z)=-\tfrac{1}{2}$, 
which is the bisector of the vertical strip $-(d-1) \leq \Re(z) \leq d-2$. 

We consider the other two classes of finite graphs. 
Let $G$ be a complete bipartite graph of type $(2,d-2)$, 
i.e., the edges of $G$ are either $\{1,j\}$ or $\{2,j\}$ with $3 \leq j \leq d.$ 
Then the $\delta$-polynomial of $\Pc_G^{\pm}$ coincides with 
$$
(1+t)^{d-3}(1+2(d-2)t+t^2). 
$$
By computational experiences, we conjecture that 
the real parts of all the roots of $i(\Pc_G^{\pm}, m)$ are equal to $-\tfrac{1}{2}$. 

Let $G$ be a complete graph with $d$ vertices and 
$\delta(\Pc_G^{\pm})=(\delta_0,\delta_1,\ldots,\delta_{d-1})$ be 
its $\delta$-vector. 
In \cite[Theorem 13]{ABHPS2008}, the $\delta(\Pc_G^{\pm})$ is calculated; that is, 
$$
\delta_i=\binom{d-1}{i}^2, i=0,1,\ldots,d-1. 
$$
By computational experiences, we also conjecture that 
the real parts of all the roots of $i(\Pc_G^{\pm}, m)$ are equal to $-\tfrac{1}{2}$. 

In addition, if $d \leq 6$, then the real parts of all the roots of $i(\Pc_G^{\pm}, m)$ 
are equal to $-\tfrac{1}{2}$ for any graph with $d$ vertices. 
However, it is not true for $d=7$ or $d=8$. In fact, there are some counterexamples. 
The following Figures \ref{fig:sep7} and \ref{fig:sep8} illustrate 
how the roots are distanced from the line $\Re(z)=-\frac{1}{2}$. 
(They are computed by {\tt CoCoA}~\cite{CocoaSystem} and {\tt Maple}~\cite{maple}.) 
\\
\begin{figure}[phtb]
\centering%
\includegraphics[scale=0.4,trim=0 0 0 300]{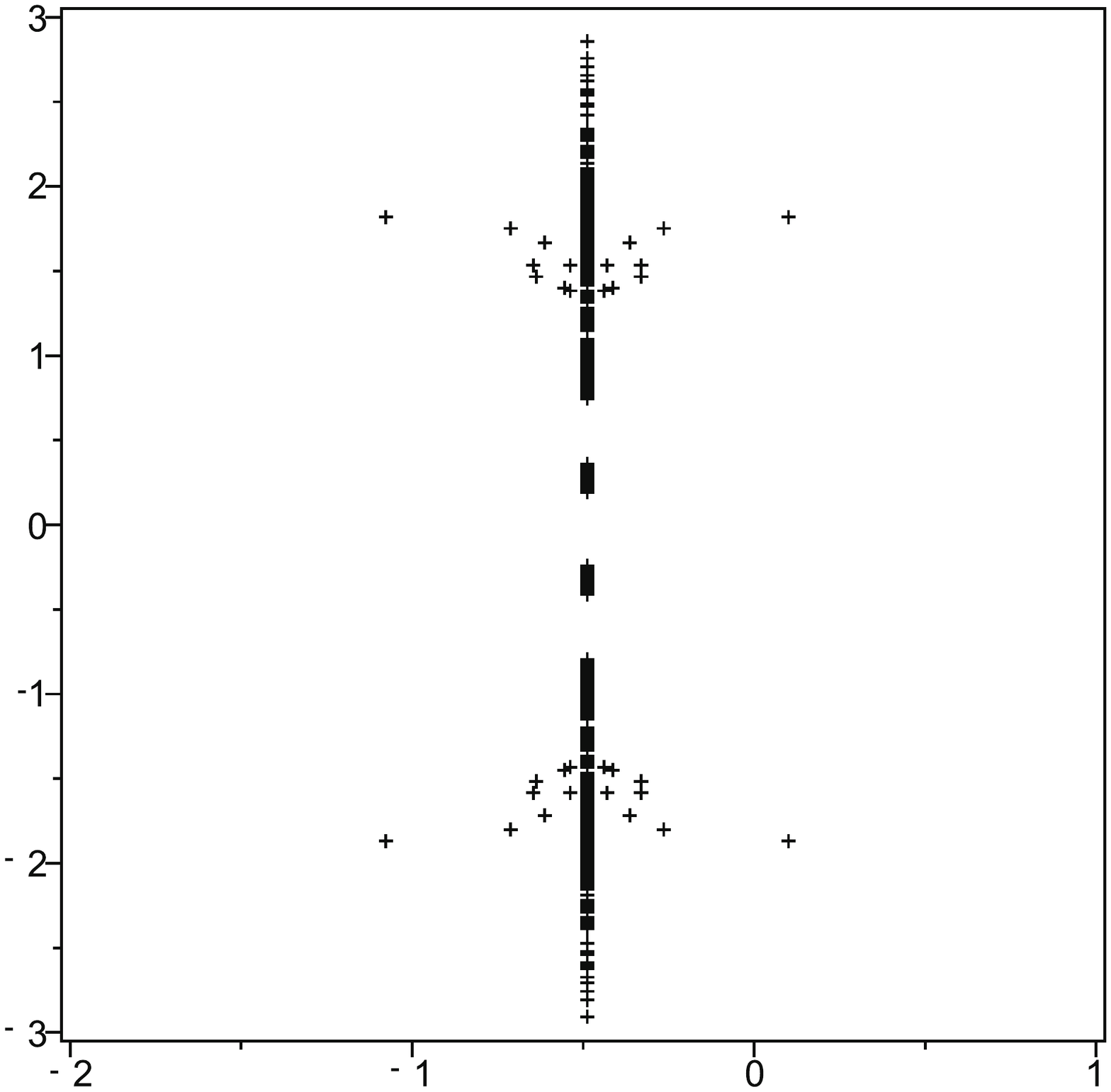}
\caption{$d=7$}
\label{fig:sep7}
\end{figure}
\begin{figure}[htbp]
\centering%
\includegraphics[scale=0.4,trim=0 0 0 250]{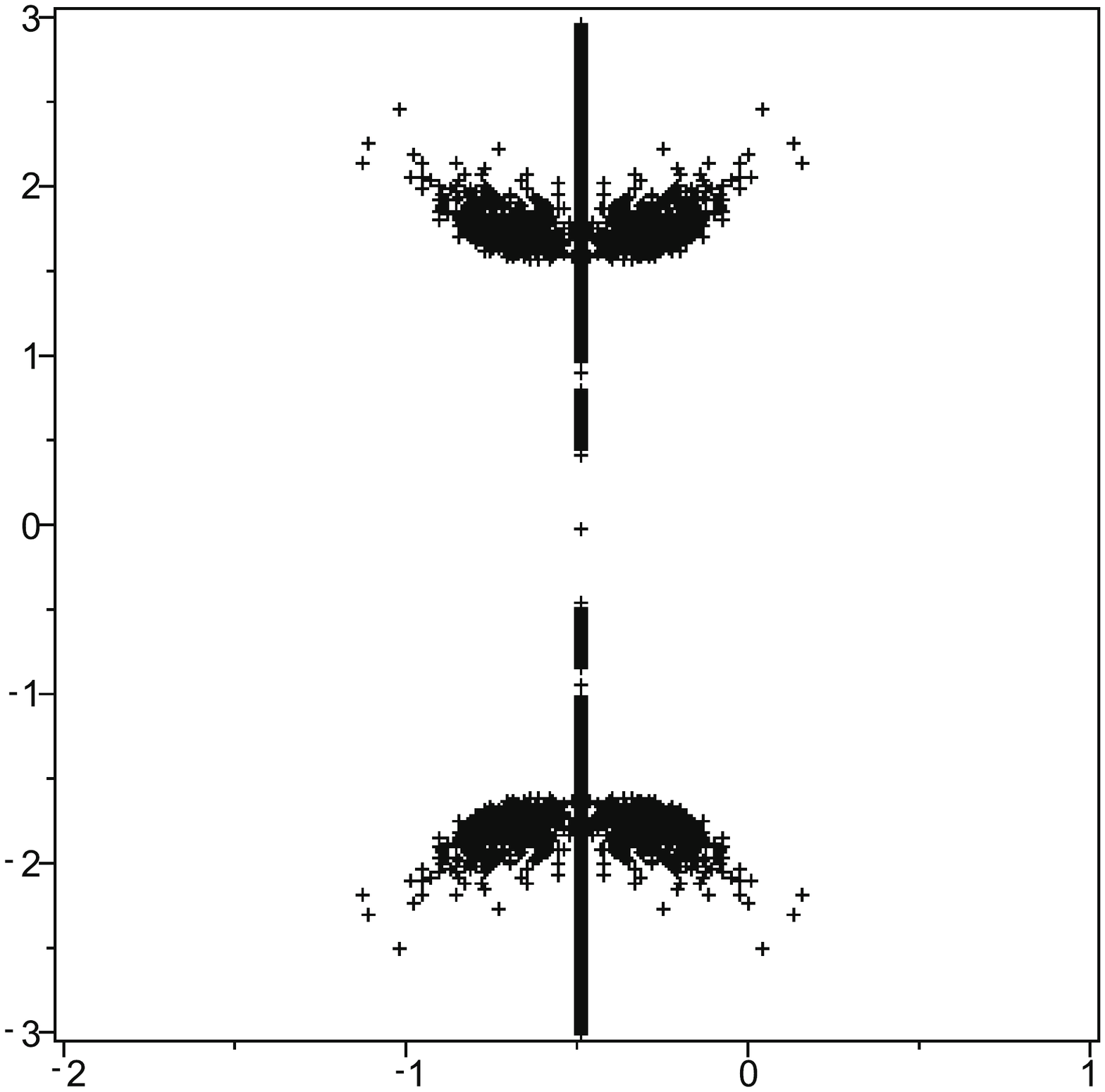}
\caption{$d=8$}
\label{fig:sep8}
\end{figure}

Let $G$ be a cycle of length $d$. 
When $d \leq 6$, although the real parts of all the roots of $i(\Pc_G^{\pm},m)$ 
are equal to $-\tfrac{1}{2}$, 
there are also some counterexamples when $d \geq 7$. 
The following Figure \ref{fig:cycle30} illustrates 
the behavior of the roots for $7 \leq d \leq 30$. 
\\
\begin{figure}[htbp]
\centering%
\includegraphics[scale=0.4,trim=0 0 0 300]{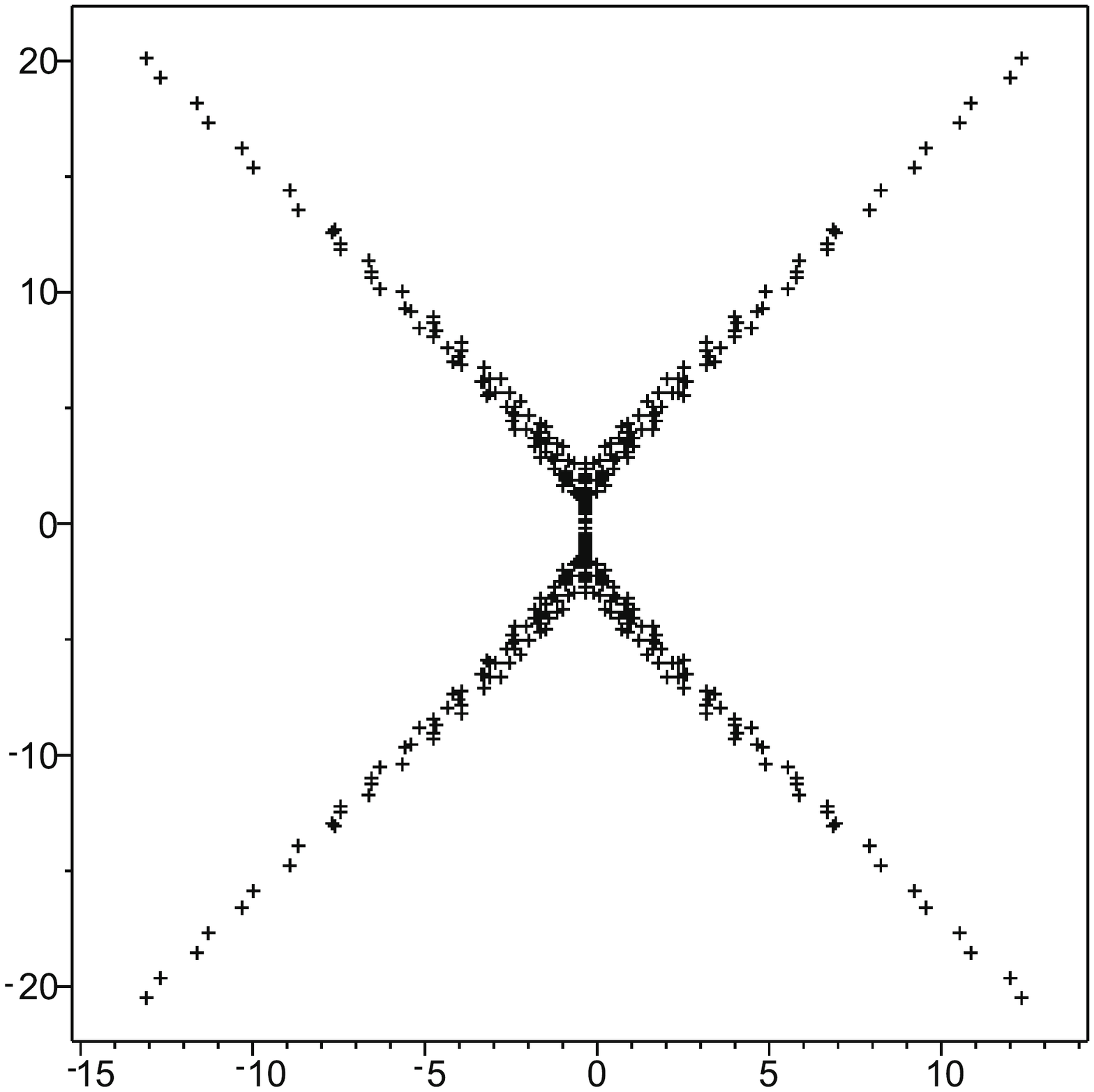}
\caption{all cycles $7 \le d \le 30$}
\label{fig:cycle30}
\end{figure}

However, all the roots $z$ of $i(\Pc_G^{\pm},m)$ whose real parts 
are not equal to $-\tfrac{1}{2}$ satisfy $-(d-1) \leq \Re(z) \leq d-2$. 
In more detail, they satisfy $-\tfrac{d-1}{2} \leq \Re(z) \leq \tfrac{d-1}{2}-1$. 
Then we propose the following: 
\begin{Conjecture}\label{conj:narrowstrip}
All roots $\a$ of the Ehrhart polynomials of 
Gorenstein Fano polytopes of dimension $D$ 
satisfy $-\tfrac{D}{2} \leq \Re(\a) \leq \tfrac{D}{2}-1$. 
\end{Conjecture}

In the table drawn below, in the second row, 
the number of connected simple graphs 
with $d ( \leq 8)$ vertices, up to isomorphism, is written. 
In the third row, among these, the number of graphs, 
up to unimodular equivalence, i.e., 
satisfying the condition in Theorem \ref{nonequivalent}, is written. 
In the fourth row, among these, in turn, the number of graphs 
that are counterexamples, i.e., 
there is a root of $i(\Pc_G^{\pm},m)$ 
whose real part is not equal to $-\tfrac{1}{2}$, is written. 

\begin{center}
\begin{tabular}{|c|c|c|c|c|c|c|c|} \hline
                 &$d=2$ &$d=3$ &$d=4$ &$d=5$ &$d=6$ &$d=7$ &$d=8$   \\ \hline
Connected graphs &$1$   &$2$   &$6$   &$21$  &$112$ &$853$ &$11117$ \\ \hline
Non equivalent   &$1$   &$2$   &$5$   &$16$  &$75$  &$560$ &$7772$  \\ \hline
Counterexamples  &$0$   &$0$   &$0$   &$0$   &$0$   &$12$  &$1092$  \\ \hline
\end{tabular}
\end{center}

\appendix
\section{Method of Computation}
\label{sec:metho}

This appendix presents an outline of the procedure used to compute the roots of 
the Ehrhart polynomials of edge or symmetric edge polytopes 
in Sections~\ref{sec:simple} and~\ref{sec:symmetric}.
Both polytopes are constructed from connected simple graphs.
For each number of vertices $d$, steps below are taken.

\begin{enumerate}
\item Construct the set of connected simple graphs of order $d$.
\item Obtain a facet representation of a polytope for each graph.
\item Compute the Hilbert series for a facet representation.
\item Build the Ehrhart polynomial from the series and solve it.
\end{enumerate}

The program for step 1 was written by the authors in the {\tt Python} programming
language with an aid of {\tt NZMATH}~\cite{NZMATH,nzmath2006}.
The source code is available at:
\begin{center}
\url{https://bitbucket.org/mft/csg/} .
\end{center}
Step 2 is performed with {\tt Polymake}~\cite{polymake,polymakewiki}.
Then, {\tt LattE}~\cite{LattE} (or {\tt LattE macchiato}~\cite{lattemacchiato})
computes the series for step 3.
The final step uses {\tt Maxima}~\cite{maxima} or
{\tt Maple}~\cite{maple}.

A small remark has to be made on the interface between steps 3 and 4. 
If one uses {\tt LattE}'s rational function as the input to {\tt Maxima},
memory consumption becomes very high. 
{\tt LattE} can send it to {\tt Maple} by itself if you
specify ``simplify,'' but this still presents the same problem for the user.
Instead, it is preferable to use the coefficient of the first several terms 
of the Taylor expansion for interpolation.

\bibliographystyle{habbrv}
\bibliography{ehrhart,graph,softwares}

\begin{thebibliography}{10}

\bibitem{ABHPS2008}
F.~Ardila, M.~Beck, S.~Ho\c{s}ten, J.~Pfeifle, and K.~Seashore.
\newblock Root polytopes and growth series of root lattices.
\newblock ar{X}iv, math/0809.5123.

\bibitem{Bat1994}
V.~Batyrev.
\newblock Dual polyhedra and mirror symmetry for {C}alabi--{Y}au hypersurfaces
  in toric varieties.
\newblock {\em J. Algebraic Geom.}, 3:493--535, 1994.

\bibitem{BDDPS2005}
M.~Beck, J.~A. De~Loera, M.~Develin, J.~Pfeifle, and R.~P. Stanley.
\newblock Coefficients and roots of {E}hrhart polynomials.
\newblock {\em Cont. Math.}, 374:15--36, 2005, math/0402148.

\bibitem{BHW2007}
C.~Bey, M.~Henk, and J.~M. Wills.
\newblock Notes on the roots of {E}hrhart polynomials.
\newblock {\em Discrete Comput. Geom.}, 38:81--98, 2007.

\bibitem{Bra2008}
B.~Braun.
\newblock Norm bounds for {E}hrhart polynomial roots.
\newblock {\em Discrete Comput. Geom.}, 39:191--193, 2008.

\bibitem{BD2006}
B.~Braun and M.~Develin.
\newblock {E}hrhart polynomial roots and {S}tanley's non-negativity theorem.
\newblock ar{X}iv, math/0610399.

\bibitem{CocoaSystem}
{CoCoA}Team.
\newblock {{\hbox{\rm C\kern-.13em o\kern-.07em C\kern-.13em o\kern-.15em A}}}:
  a system for doing {C}omputations in {C}ommutative {A}lgebra.
\newblock \url{http://cocoa.dima.unige.it/}.

\bibitem{LattE}
J.~A. De~Loera, D.~Haws, R.~Hemmecke, P.~Huggins, J.~Tauzer, and R.~Yoshida.
\newblock {\tt {L}att{E}}.
\newblock \url{http://wwww.math.ucdavis.edu/~latte/}.

\bibitem{polymake}
E.~Gawrilow and M.~Joswig.
\newblock polymake: a framework for analyzing convex polytopes.
\newblock In G.~Kalai and G.~M. Ziegler, editors, {\em Polytopes ---
  Combinatorics and Computation}, pages 43--74. Birkh\"auser, 2000.

\bibitem{Harary}
F.~Harary.
\newblock {\em Graph Theory}.
\newblock Addison-Wesley, Reading, 1969.

\bibitem{HarPal}
F.~Harary and E.~M. Palmer.
\newblock {\em Graphical Enumeration}.
\newblock Academic Press, New York and London, 1973.

\bibitem{HSW2005}
M.~Henk, A.~Sch\"urmann, and J.~M. Wills.
\newblock {E}hrhart polynomials and successive minima.
\newblock {\em Mathematika}, 52:1--16, 2005.

\bibitem{Hibi1992}
T.~Hibi.
\newblock {\em Algebraic combinatorics of convex polytopes}.
\newblock Carslaw Publications, Glebe, N.S.W., 1992.

\bibitem{Hib1992}
T.~Hibi.
\newblock Dual polytopes of rational convex polytopes.
\newblock {\em Combinatorica}, 12:237--240, 1992.

\bibitem{Hib1994}
T.~Hibi.
\newblock A lower bound theorem for {E}hrhart polynomials of convex polytopes.
\newblock {\em Adv. in Math.}, 105:162--165, 1994.

\bibitem{lattemacchiato}
M.~K\"{o}ppe.
\newblock {\tt {L}att{E} macchiato}.
\newblock \url{http://www.math.ucdavis.edu/~mkoeppe/latte/}.

\bibitem{nzmath2006}
T.~Matsui.
\newblock Development of {N}{Z}{M}{A}{T}{H}.
\newblock In A.~Iglesias and N.~Takayama, editors, {\em Mathematical Software
  --- ICMS 2006}, volume 4151 of {\em Lecture Notes in Computer Science}, pages
  158--169. Springer-Verlag, September 2006.

\bibitem{maxima}
Maxima.sourceforge.net.
\newblock {\tt Maxima}, a computer algebra system.
\newblock \url{http://maxima.sourceforge.net/}.

\bibitem{NZMATH}
{NZMATH development group}.
\newblock {\tt {N}{Z}{M}{A}{T}{H}}.
\newblock \url{http://tnt.math.metro-u.ac.jp/nzmath/}.

\bibitem{OH1998}
H.~Ohsugi and T.~Hibi.
\newblock Normal polytopes arising from finite graphs.
\newblock {\em J. Algebra}, 207:209--426, 1998.

\bibitem{OH2000}
H.~Ohsugi and T.~Hibi.
\newblock Compressed polytopes, initial ideals and complete multipartite
  graphs.
\newblock {\em Illinois Journal of Mathematics}, 44(2):391--406, 2000.

\bibitem{OH2008}
H.~Ohsugi and T.~Hibi.
\newblock Simple polytopes arising from finite graphs.
\newblock In {\em ITSL}, pages 73--79, 2008.
\newblock also available at ar{X}iv:math/0804.4287.

\bibitem{Pfe2010}
J.~Pfeifle.
\newblock Gale duality bounds for roots of polynomials with nonnegative
  coefficients.
\newblock {\em J. Combin. Theory Ser. A}, 117(3):248--271, 2010.

\bibitem{polymakewiki}
{\tt polymake}.
\newblock \url{http://www.opt.tu-darmstadt.de/polymake/doku.php}.

\bibitem{Rod2002}
F.~Rodriguez-Villegas.
\newblock On the zeros of certain polynomials.
\newblock {\em Proc. Amer. Math. Soc.}, 130:2251--2254, 2002.

\bibitem{Schrijver1986}
A.~Schrijver.
\newblock {\em Theory of Linear and Integer Programming}.
\newblock John Wiley \& Sons, 1986.

\bibitem{Sta1980}
R.~P. Stanley.
\newblock Decompositions of rational convex polytopes.
\newblock {\em Annals of Discrete Math.}, 6:333--342, 1980.

\bibitem{Stanley1986}
R.~P. Stanley.
\newblock {\em Enumerative Combinatorics, volume I}.
\newblock Wadsworth \& Brooks / Cole Advanced Books, Monterey, Calif., 1986.

\bibitem{Sta1993}
R.~P. Stanley.
\newblock A monotonicity property of $h$-vectors and $h^*$-vectors.
\newblock {\em Europ. J. Combinatorics}, 14:251--258, 1993.

\bibitem{Sturmfels1995}
B.~Sturmfels.
\newblock {\em Gr\"obner Bases and Convex Polytopes}, volume~8 of {\em
  University Lecture Series}.
\newblock American Mathematical Society, Providence, RI, 1995.
\newblock ISBN 0-8218-0487-1.

\bibitem{maple}
{Waterloo Maple Inc.}
\newblock {\tt Maple}.
\newblock \url{http://www.maplesoft.com/products/Maple/}.

\bibitem{Wilson}
R.~J. Wilson.
\newblock {\em Introduction to Graph Theory}.
\newblock Addison-Wesley, Reading, fourth edition, 1996.

\end{thebibliography}

\end{document}